# DIFFUSION APPROXIMATIONS FOR CONTROLLED STOCHASTIC NETWORKS: AN ASYMPTOTIC BOUND FOR THE VALUE FUNCTION

By Amarjit Budhiraja[1] and Arka Prasanna Ghosh[1]

*University of North Carolina at Chapel Hill and Iowa State University*

We consider the scheduling control problem for a family of unitary networks under heavy traffic, with general interarrival and service times, probabilistic routing and infinite horizon discounted linear holding cost. A natural nonanticipativity condition for admissibility of control policies is introduced. The condition is seen to hold for a broad class of problems. Using this formulation of admissible controls and a time-transformation technique, we establish that the infimum of the cost for the network control problem over all admissible sequencing control policies is asymptotically bounded below by the value function of an associated diffusion control problem (the Brownian control problem). This result provides a useful bound on the best achievable performance for any admissible control policy for a wide class of networks.

**1. Introduction.** In [12], Harrison introduced "controlled Brownian networks" as *formal* approximations to stochastic network models under diffusion scaling. Since then, several authors have used such Brownian network models and the corresponding control problem [the so-called Brownian control problem (BCP)] in obtaining control policies for the underlying queueing networks. In particular, there have been several works in recent years [1, 2, 5, 7, 20] that consider specific network models for which the corresponding BCP is explicitly solvable and, guided by the solution of the BCP, propose a control policy for the network which is then shown to be asymptotically optimal. However, there is a critical lack of mathematical theory which provides a rigorous basis for using the BCP as an approximating model for

Received July 2005; revised April 2006.
[1]Supported in part by ARO Grant W911NF-04-1-0230.

a general family of controlled queueing networks. In this paper, we will consider a broad family of unitary networks (cf. [5]) under heavy traffic with general interarrival, service times and probabilistic routing with an infinite horizon discounted linear holding cost. We will prove that if a control policy is *admissible* (see Definition 2.6) for the network control problem, then the associated cost of using this policy is asymptotically bounded below by the value function of the corresponding diffusion control problem. This is the main result of the paper and is given in Theorem 3.1 (and in Corollary 3.2). Although the result is far from optimal—in particular one would also like to establish the reverse inequality—it is the first result, in the broad generality we consider, that establishes any asymptotic relationship between the controlled stochastic network and its formal diffusion approximation. Indeed, even a proper formulation of admissible control policies for queueing networks that allows for weak convergence analysis to asymptotically relate it to the BCP has been critically lacking. Such a formulation, which in particular needs to include appropriate nonanticipativity conditions, is introduced in Definition 2.6. It will be seen in Theorem 5.4 that there exists a very large class of natural control policies which are admissible in the sense of Definition 2.6. (See also Proposition 2.8 for another sufficient condition for admissibility.) Roughly speaking, a control policy that does not change between two successive event times (an *event* is defined to be an exogenous arrival into the system or a completion of a service by some server), and at any event instant is a function of all the previous interarrival times, service times and routing decisions, is seen to satisfy the key requirements of admissibility in Definition 2.6. Once a proper formulation of admissibility is available, it becomes possible to use the machinery of weak convergence theory to relate the asymptotic value of the queuing control problem to the value function of the associated formal diffusion approximation. The main obstacle, in making rigorous the formal heavy traffic limit needed in the diffusion approximation approach, is proving the tightness of the control terms in the state evolution equation; see, for example, (2.22). Although for each $r$ the control term $\hat{Y}^r$ on the right-hand side of (2.22) is Lipschitz continuous, its Lipschitz constant increases to $\infty$, linearly in $r$ as $r \to \infty$. One indication that the asymptotic analysis is somewhat delicate is that although the control terms for the network control problem are continuous (in fact Lipschitz) for each fixed $r$, an optimal control in the limiting Brownian control problem can in fact have jumps. In this paper, we develop a general approach to argue the convergence of the cost function which does not rely on the tightness of the various processes describing the dynamics, but rather on the tightness of suitable time-rescaled versions of the same. Although such time rescaling ideas go back to Meyer and Zheng [19] and Kurtz [15], their application in stochastic control problems was pioneered by Kushner and coauthors [17, 18]. Indeed, inspiration for many ideas in the current paper comes from these latter works. In a related paper, Kurtz and Stockbridge [16] consider a control problem for a stable Jackson network with Markovian primitives. They are concerned with an ergodic cost problem associated with the approximating controlled diffusion; however, their techniques and goals are quite different from those of the current paper.

The paper is organized as follows. The next section describes the queueing network and the associated control problem. Our key nonanticipativity condition on control policies is introduced in Definition 2.6. Two diffusion control problems, BCP and the corresponding equivalent workload formulation (EWF), which arise from formal diffusion approximations of the network control problem are introduced as well. These control problems are given a "weak formulation" which is natural in view of the weak convergence ideas underpinning all heavy traffic approximation methods. Section 3 contains the main result of this paper (Theorem 3.1). The proofs of all but one of the key results, namely Theorem 3.7, are contained in Section 3 as well. This latter theorem is at the heart of our analysis and its proof is given separately in Section 4. The main ingredients of the proof are the time transformation in Lemma 3.6 (which allows for appropriate tightness estimates) and relationships between various multiparameter filtrations, stopping times and martingales associated with the queuing model. A careful analysis of the latter enables one to show that a typical limit point of the (time-transformed) control sequence has suitable adaptability properties, and the limit points of the "pre-Brownian motions" in the queuing model are martingales with respect to the desired filtration. We refer the reader to Remark 3.8 for more on these key adaptability issues. Finally, in Section 5 of this paper, we show that there exists a broad family of control policies which satisfy the admissibility requirements of Definition 2.6.

The following notation will be used. The space of reals (nonnegative reals), positive integers (nonnegative integers) will be denoted by $\mathbb{R}$ ($\mathbb{R}_+$), $\mathbb{N}$ ($\mathbb{N}_0$), respectively. For $m \geq 1$, $\mathcal{C}^m$ will denote the space of continuous functions from $[0, \infty)$ to $\mathbb{R}^m$. For $m \geq 1$, $\mathcal{D}^m$ will denote the space of right continuous functions with left limits, from $[0, \infty)$ to $\mathbb{R}^m$ with the usual Skorohod topology and $\mathcal{B}(\mathcal{D}^m)$ the corresponding Borel sigma-field. All continuous-time processes considered in this paper will have sample paths in $\mathcal{D}^m$. All vectors will be column vectors and all vector inequalities are to be interpreted componentwise. We will call a function $f \in \mathcal{D}^m$ *nonnegative* if $f(t) \geq 0$ for all $t \in \mathbb{R}_+$. A function $f \in \mathcal{D}^m$ is called *nondecreasing* if it is nondecreasing in each component. The weak convergence of processes $Z_n$ to $Z$ as elements of $\mathcal{D}^m$ will be denoted by $Z_n \Rightarrow Z$. A sequence of processes $\{Z_n\}$ is tight if and only if the measures induced by the $Z_n$'s on $(\mathcal{D}^m, \mathcal{B}(\mathcal{D}^m))$ form a tight sequence. A sequence of processes with paths in $\mathcal{D}^m$ ($m \geq 1$) is called $\mathcal{C}$-*tight* if it is tight in $\mathcal{D}^m$ and any weak limit point of the sequence has paths in $\mathcal{C}^m$ almost surely. For processes $\{Z_n\}$, $Z$ defined on a common probability space, we say that $Z_n$ converge to $Z$ uniformly on compact time intervals (u.o.c.) in probability (almost surely) if for all $t > 0$, $\sup_{0 \leq s \leq t} |Z_n(s) - Z(s)|$ converges to zero in probability a.s. The inherent nature of the problem makes this article heavy in notation. To ease the notational burden, standard notation (see, e.g., [6]) for different processes is used (e.g., $Q$ for queue-length, $I$ for idle time, $W$ for workload process, etc.). We also use standard notation, for example, $\bar{W}, \hat{W}$, to denote fluid scaled, respectively, diffusion scaled, versions of various processes of interest. Time rescaling ideas used in the paper lead to processes [see (3.36)] that are obtained via time transformation of the original processes of interest; these are typically denoted with a superscript, as in $\check{W}$. Since we deal with vector-valued processes, we use both $(a_i, i = 1, \ldots, m; b_i, i = 1, \ldots, n)$ and $(a_1, \ldots, a_m, b_1, \ldots, b_n)$ to denote a vector the first $m$ components of which are given by $a_1, \ldots, a_m$ and the next $n$ components of which are given by $b_1, \ldots, b_n$, and so on. The notation $(c_i^j, i = 1, \ldots, m, j = 1, \ldots, n)$ corresponds to the vector $(c_1^1, \ldots, c_m^1, c_1^2, \ldots, c_m^2, \ldots, c_1^n, \ldots, c_m^n)$.

**2. Stochastic network model and the control problem.** Let $(\Omega, \mathcal{F}, \mathbb{P})$ be a probability space. All the random variables associated with the network model described below are assumed to be defined on this probability space, unless specified otherwise. The expectation operation under $\mathbb{P}$ will be denoted by $\mathbb{E}$.

*Network structure.* We begin by introducing the family of queuing network models that will be considered in this work. Figure 1 gives a schematic for such a model. The network has $\mathbf{I}$ infinite capacity buffers (to store $\mathbf{I}$ different classes of jobs) and $\mathbf{K}$ nonidentical servers for processing jobs. Arrivals of jobs can be from outside the system and/or from the internal rerouting of jobs that have already been served by some server. Several different servers may process jobs from a particular buffer. Service from a given buffer $i$ by a given server $k$ is called an activity. It is assumed that there are $\mathbf{J}$ activities [at most one activity for a server–buffer pair $(i, k)$, so that $\mathbf{J} \leq \mathbf{I} \cdot \mathbf{K}$]. The activities are labeled $j = 1, \ldots, \mathbf{J}$. It is also assumed that the integers $\mathbf{I}, \mathbf{J}, \mathbf{K}$ are strictly positive.

The correspondence between activities and buffers, and activities and servers, are described by two matrices, $C$ and $A$, respectively. $C$ is an $\mathbf{I} \times \mathbf{J}$ matrix with $C_{ij} = 1$ if the $j$th activity processes jobs from buffer $i$, and $C_{ij} = 0$ otherwise. $A$ is a $\mathbf{K} \times \mathbf{J}$ matrix with $A_{kj} = 1$ if the $k$th server is associated with the $j$th activity, and $A_{kj} = 0$ otherwise. Each activity associates one buffer and one server, and so each column of $C$ has exactly one 1 (and similarly, each column of $A$ has exactly one 1). We will further assume that each row of $C$ (resp. $A$) has at least one 1, that is, each buffer is processed by (resp. server is processing) at least one activity.

It is also assumed that once a job starts being processed by an activity, it must complete its service with that activity, even if its service is interrupted for some time (e.g., for preemption by a job from another buffer). When service of a partially completed job is resumed, it is resumed from the point of preemption—that is the job needs only the remaining service time from the server to be completed (preemptive-resume policy). Also, an activity must complete the service of any job that it has started before starting another job from the same buffer. An activity always selects the oldest job in the buffer that is not yet served when starting a new service [i.e., First In First Out (FIFO) within class].

*Stochastic primitives.* The networks considered in this paper will be nearly critically loaded, that is, "approaching heavy traffic," as made precise in Assumption 2.3 below. Mathematically, this is modeled by considering a sequence of networks $\{\mathcal{N}^r\}$, each network in the sequence having identical structure except for the rate parameters that may depend on $r$. Here, $r \in \mathbb{S} \subseteq \mathbb{R}^+$, where $\mathbb{S}$ is a countable set: $\{r_1, r_2, \ldots\}$ with $1 \leq r_1 < r_2 < \cdots$ and $r_n \to \infty$, as $n \to \infty$. For notational simplicity, throughout the paper, we will write the limit along the sequence $r_n$ as $n \to \infty$ simply as "$r \to \infty$." Also, $r$ will always be taken to be an element of $\mathbb{S}$ and thus hereafter the qualifier $r \in \mathbb{S}$ will not be stated explicitly. The parameters of the networks are assumed to satisfy the heavy traffic condition (see Definition 2.2). The physical network that is "close to critical loading" corresponds to one of the networks in the sequence, for a large value of $r$.

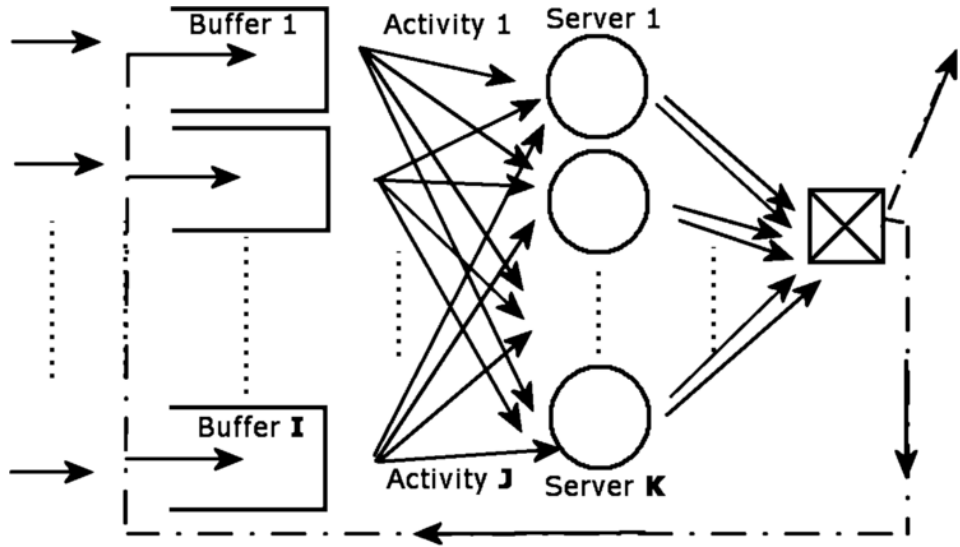

Fig. 1. *A unitary network.*

The $r$th network $\mathcal{N}^r$ is described as follows. If the $i$th class of jobs ($i = 1, \ldots, \mathbf{I}$) has exogenous arrivals, the interarrival times are given by a sequence of nonnegative random variables $\{u_i^r(n) : n \geq 1\}$ that are i.i.d. with mean and standard deviation $1/\alpha_i^r$ and $\sigma_i^{u,r} \in (0, \infty)$, respectively. Let, by relabeling if necessary, the buffers with exogenous arrivals correspond to $i = 1, \ldots, \mathbf{I}'$, where $\mathbf{I}' \leq \mathbf{I}$. We set $\alpha_i^r, \sigma_i^{u,r} = 0$ and $u_i^r(n) = \infty$, $n \geq 1$, for $i = \mathbf{I}' + 1, \ldots, \mathbf{I}$. Service times of the $j$th type of activity (for $j = 1, \ldots, \mathbf{J}$) are given by a sequence of nonnegative random variables $\{v_j^r(n) : n \geq 1\}$ that are i.i.d. with mean and standard deviation $1/\beta_j^r$ and $\sigma_j^{v,r} \in (0, \infty)$, respectively. We will assume that the above random variables are in fact strictly positive, that is,

(2.1)  For all $i = 1, \ldots, \mathbf{I}$, $j = 1, \ldots, \mathbf{J}$   $\mathbb{P}(u_i^r(1) > 0) = \mathbb{P}(v_j^r(1) > 0) = 1$.

We will further impose the following uniform integrability condition: For all $j = 1, \ldots, \mathbf{J}$ and all $i = 1, \ldots, \mathbf{I}'$,

(2.2) the sequences $\{(u_i^r(1))^2\}_r$ and $\{(v_j^r(1))^2\}_r$ are uniformly integrable.

Rerouting of jobs completed by the $j$th activity is specified by a sequence of $(\mathbf{I} + 1)$-dimensional vectors $\{(\phi_0^{j,r}(n), \phi^{j,r}(n)), n \geq 1\}$, where $\phi^{j,r}(n) = (\phi_i^{j,r}(n) : i = 1, \ldots, \mathbf{I})$. For each $j$ and $i = 0, 1, \ldots, \mathbf{I}$, $\phi_i^{j,r}(n) = 1$ if the $n$th completed job by activity $j$ gets rerouted to buffer $i$, and takes the value zero otherwise, where $i = 0$ represents jobs leaving the system. It is assumed that for each fixed $r$, $\{(\phi_0^{j,r}(n), \phi^{j,r}(n)), n \geq 1\}$, for $j = 1, \ldots, \mathbf{J}$ are independent sequences of i.i.d. $Multinomial_{(\mathbf{I}+1)}(1, (p_0^j, p^j))$, where $p^j = (p_i^j : i = 1, \ldots, \mathbf{J})$. That, in particular, means for $j = 1, \ldots, \mathbf{J}$, $n \geq 1$, $\sum_{i=0}^{\mathbf{I}} \phi_i^{j,r}(n) = \sum_{i=0}^{\mathbf{I}} p_i^j = 1$. Furthermore, for fixed $j = 1, \ldots, \mathbf{J}, i_1, i_2 = 1, \ldots, \mathbf{I}$, if $\sigma_{i_1 i_2}^{\phi_j} \doteq \mathrm{Cov}(\phi_{i_1}^{j,r}(n), \phi_{i_2}^{j,r}(n))$, then

(2.3) $$\sigma_{i_1 i_2}^{\phi_j} = \begin{cases} p_i^j(1 - p_i^j), & \text{if } i_1 = i_2 = i, \\ -p_{i_1}^j p_{i_2}^j, & \text{if } i_1 \neq i_2. \end{cases}$$

We also assume that for $i = 1, \ldots, \mathbf{I}$, $j = 1, \ldots, \mathbf{J}$, the sequences $\{u_i^r(n) : n \geq 1\}$, $\{v_j^r(n) : n \geq 1\}$ and $\{\phi^{j,r}(n), n \geq 1\}$ are mutually independent.

Next, we introduce the primitive stochastic processes $(E^r, S^r)$ that describe the state dynamics for the model. The process $(E_1^r, \ldots, E_{\mathbf{I}'}^r)$ is the $\mathbf{I}'$-dimensional exogenous arrival process, that is, for each $i = 1, \ldots, \mathbf{I}'$, $E_i^r(t)$ is a renewal process which denotes the number of jobs that have arrived to buffer $i$ from outside the system in $[0, t]$. For a class $i$ to which there are no exogenous arrivals (i.e., $i = \mathbf{I}' + 1, \ldots, \mathbf{I}$), we set $E_i^r(t) = 0$ for all $t \geq 0$. We will denote the process $(E_1^r, \ldots, E_{\mathbf{I}}^r)'$ by $E^r$. For each activity $j = 1, \ldots, \mathbf{J}$, there is a renewal process $S_j^r$ such that $S_j^r(t)$ denotes the number of complete jobs that could be processed by activity $j$ in $[0,t]$ *if the associated server worked continuously and exclusively on jobs from the associated buffer in* $[0,t]$. The vector $(S_1^r, \ldots, S_\mathbf{J}^r)'$ is denoted by $S^r$. More precisely, for $i = 1, \ldots, \mathbf{I}$, $j = 1, \ldots, \mathbf{J}$, $m_i \geq 1, n_j \geq 1$, let

$$\xi_i^r(m_i) \doteq \sum_{n=1}^{m_i} u_i^r(n), \qquad \eta_j^r(n_j) \doteq \sum_{n=1}^{n_j} v_j^r(n). \tag{2.4}$$

We set $\xi_i^r(0) = 0$ and $\eta_j^r(0) = 0$. Then $E_i^r$, $S_j^r$ are renewal processes given as follows:

$$\begin{aligned} E_i^r(t) &= \max\{m_i \geq 0 : \xi_i^r(m_i) \leq t\}, \\ S_j^r(t) &= \max\{n_j \geq 1 : \eta_j^r(n_j) \leq t\}, \qquad t \geq 0. \end{aligned} \tag{2.5}$$

Finally, we introduce the routing sequences. Let $\boldsymbol{\Phi}_i^{j,r}(n)$ denote the number of jobs that are routed to the $i$th buffer, out of the first $n$ jobs completed by activity $j$. Thus, for $i = 1, \ldots, \mathbf{I}$, $j = 1, \ldots, \mathbf{J}$,

$$\boldsymbol{\Phi}_i^{j,r}(n) = \sum_{m=1}^{n} \phi_i^{j,r}(m), \qquad n = 1, 2, \ldots. \tag{2.6}$$

We will denote by $\{\boldsymbol{\Phi}^{j,r}(n)\}$ the $\mathbf{I}$-dimensional sequence $\{(\boldsymbol{\Phi}_1^{j,r}(n), \ldots, \boldsymbol{\Phi}_\mathbf{I}^{j,r}(n))'\}$ corresponding to the routing of jobs completed by the $j$th activity. Also, the matrix $(\boldsymbol{\Phi}^{1,r}(n), \boldsymbol{\Phi}^{2,r}(n), \ldots, \boldsymbol{\Phi}^{\mathbf{J},r}(n))$ will be denoted by $\boldsymbol{\Phi}^r(n)$.

*Control.* A *scheduling policy or control* is specified by a nonnegative, nondecreasing $\mathbf{J}$-dimensional process $T^r = \{(T_1^r(t), \ldots, T_\mathbf{J}^r(t)), \ t \geq 0\}$. For any $j = 1, \ldots, \mathbf{J}$, $t \geq 0$, $T_j^r(t)$ represents the cumulative amount of time spent on the $j$th activity up to time $t$. For a control $T^r$ to be admissible, it is required to satisfy additional properties which are specified below in Definition 2.6.

*State processes.* For a given scheduling policy $T^r$, the state processes of the network are the associated $\mathbf{J}$-dimensional queue length process $Q^r$ and the $\mathbf{K}$-dimensional idle time process $I^r$. For each $t \geq 0$, $i = 1, \ldots, \mathbf{I}$, $Q_i^r(t)$ represents the queue length at the $i$th buffer at time $t$ (including the jobs that are in service at that time), and for $k = 1, \ldots, \mathbf{K}$, $I_k^r(t)$ is the total amount of time that the $k$th server has idled, up to time $t$. Let $q^r = Q^r(0) \in \mathbb{R}_+^\mathbf{J}$ be the $\mathbf{J}$-dimensional vector of initial queue lengths. Note that, for $j = 1, \ldots, \mathbf{J}$, $t \geq 0$, $S_j^r(T_j^r(t))$ is the total number of services completed by the $j$th activity up to time $t$. Thus, the total number of completed jobs (by activity $j$) up to time $t$ that get rerouted to buffer $i = 1, \ldots, \mathbf{I}$ is given by the process $\boldsymbol{\Phi}_i^{j,r}(S_j^r(T_j^r(t)))$. Hence, the state of the system at time $t \geq 0$ is described by the following equations:

$$Q_i^r(t) = q^r + E_i^r(t) - \sum_{j=1}^{\mathbf{J}} C_{ij} S_j^r(T_j^r(t)) + \sum_{j=1}^{\mathbf{J}} \boldsymbol{\Phi}_i^{j,r}(S_j^r(T_j^r(t))),$$

(2.7)
$$i = 1, \ldots, \mathbf{I},$$

(2.8) $\quad I_k^r(t) = t - \sum_{j=1}^{\mathbf{J}} A_{kj} T_j^r(t), \qquad k = 1, \ldots, \mathbf{K}.$

We will only be concerned with admissible policies (see Definition 2.6) which, in particular, ensures that the queue-length and idle-time processes $(Q^r, I^r)$ are nonnegative.

*Heavy traffic.* We begin with a condition on convergence of various parameters in the sequence of networks $\{\mathcal{N}^r\}$.

ASSUMPTION 2.1. There are nonnegative $q$, $\alpha$, $\sigma^u$, $\theta_1 \in \mathbb{R}^{\mathbf{I}}$; $\beta$, $\sigma^v$, $\theta_2 \in \mathbb{R}^{\mathbf{J}}$ such that $q \geq 0$, $\alpha \geq 0$, $\beta > 0$, $\sigma^u \geq 0$, $\sigma^v > 0$ ($\alpha_i$, $\sigma_i^u = 0$ if and only if $i = \mathbf{I}' + 1, \ldots, \mathbf{I}$), and

$$\theta_1^r \doteq r(\alpha^r - \alpha) \to \theta_1, \qquad \theta_2^r \doteq r(\beta^r - \beta) \to \theta_2 \qquad \text{as } r \to \infty,$$

(2.9) $\quad \sigma^{u,r} \to \sigma^u, \qquad\qquad\qquad \sigma^{v,r} \to \sigma^v \qquad\qquad\quad \text{as } r \to \infty,$

$$\hat{q}^r \doteq \frac{q^r}{r} \to q \qquad\qquad\qquad\qquad\qquad\qquad\qquad \text{as } r \to \infty.$$

Next, we present the key heavy traffic condition, as introduced in [11] (also see [5, 6, 13]), on the sequence $\{\mathcal{N}^r\}$.

DEFINITION 2.2 (*Heavy traffic*). Define $\mathbf{I} \times \mathbf{J}$ matrices $P'$ and $R$, such that $P'_{ij} \doteq p_i^j$, for $i = 1, \ldots, \mathbf{I}$, $j = 1, \ldots, \mathbf{J}$ and

(2.10) $$R \doteq (C - P') \operatorname{diag}(\beta).$$

We say that the sequence $\{\mathcal{N}^r\}$ approaches *heavy traffic* as $r \to \infty$ if, in addition to Assumption 2.1, the following two conditions hold:

(i) There is a unique optimal solution $(x^*, \rho^*)$ to the following linear program (LP):

(2.11) $\qquad\qquad\text{minimize } \rho \quad \text{s.t.} \quad Rx = \alpha, \qquad Ax \leq \rho e, x \geq 0.$

(ii) The unique solution $(x^*, \rho^*)$ of the above linear program satisfies $\rho^* = 1$ and $Ax^* = e$.

Here $e$ is a $\mathbf{K}$-dimensional vector of ones, and for a vector $a$, $\operatorname{diag}(a)$ denotes the diagonal matrix such that the vector of its diagonal entries is $a$.

ASSUMPTION 2.3. *The sequence of networks $\{\mathcal{N}^r\}$ approaches heavy traffic as $r \to \infty$.*

Assumption 2.3 will be a standing assumption for this paper and will not be mentioned explicitly in the statements of various results.

REMARK 2.4. From Assumption 2.3, $x^*$ given in (i) of Definition 2.2 is the unique **J**-dimensional nonnegative vector satisfying

$$(2.12) \qquad Rx^* = \alpha, \qquad Ax^* = e.$$

Following [11], assume without loss of generality (by relabeling activities, if necessary), that the first **B** components of $x^*$ are strictly positive (corresponding activities are referred to as *basic*) and the rest are zero (*nonbasic* activities). For later use, we partition the following matrices and vectors in terms of *basic* and *non-basic* components:

$$(2.13) \qquad x^* = \begin{bmatrix} x_b^* \\ 0 \end{bmatrix}, \qquad T^r(\cdot) = \begin{bmatrix} T_b^r(\cdot) \\ T_n^r(\cdot) \end{bmatrix}, \qquad A = [B : N],$$

where **0** is a $(\mathbf{J} - \mathbf{B})$-dimensional vector of zeros and $B$ and $N$ are $\mathbf{K} \times \mathbf{B}$ and $\mathbf{K} \times (\mathbf{J} - \mathbf{B})$ matrices, respectively.

*Other processes.* The vector $x^*$ defined above gives the nominal allocation rates for the **J** activities. Define the following *deviation process* as the difference between $T^r$ and the nominal allocation:

$$(2.14) \qquad Y^r(t) \doteq x^* t - T^r(t), \qquad t \geq 0.$$

It follows from (2.8) and (2.12) that the idle-time process $I^r$ has the following representation:

$$I^r(t) = AY^r(t), \qquad t \geq 0.$$

Next, we define a $(\mathbf{K} + \mathbf{J} - \mathbf{B}) \times \mathbf{J}$ matrix $K$ and a $(\mathbf{K} + \mathbf{J} - \mathbf{B})$-dimensional process $U^r$ as follows:

$$(2.15) \qquad K \doteq \begin{bmatrix} B & N \\ 0 & -\mathbb{I} \end{bmatrix}, \qquad U^r(t) \doteq KY^r(t), \qquad t \geq 0,$$

where $\mathbb{I}$ above denotes a $(\mathbf{J} - \mathbf{B}) \times (\mathbf{J} - \mathbf{B})$ identity matrix. Note that

$$(2.16) \qquad U^r(t) = \begin{bmatrix} I^r(t) \\ T_n^r(t) \end{bmatrix}, \qquad t \geq 0.$$

Another process that plays a crucial role in our analysis is the workload process. This is an **L**-dimensional process ($\mathbf{L} \leq \mathbf{I}$) defined as a suitable linear combination of the queue length processes as follows:

$$(2.17) \qquad W^r(t) = \Lambda Q^r(t), \qquad t \geq 0,$$

where $\Lambda$ is an $\mathbf{L} \times \mathbf{I}$-dimensional matrix, called the *workload matrix*, whose rows are determined by the optimal solution of the dual of the linear program in (2.11) (see [13]). We will work with one of the canonical choices of the workload matrix as defined in [13]. For this choice of $\Lambda$, one can obtain (see [6]) a *nonnegative* $\mathbf{L} \times (\mathbf{K} + \mathbf{J} - \mathbf{B})$ matrix $G$ such that

$$\Lambda R = GK. \tag{2.18}$$

We make the following additional assumption on $G$:

ASSUMPTION 2.5. Each column of $G$ has at least one strictly positive entry: equivalently, there exists a $c > 0$ such that for every $u \in \mathbb{R}_+^{(\mathbf{K}+\mathbf{J}-\mathbf{B})}$,

$$|Gu| \geq c|u|.$$

The above will also be a standing assumption, so explicit reference to it will be omitted.

*Rescaled processes.* The two scalings which are exploited in our analysis are the fluid scaling (corresponding to the law of large numbers) and diffusion scaling (corresponding to the central limit theorem). These are described as follows:

*Fluid scaled processes.* These are obtained from the original process by accelerating time by a factor of $r^2$ and dividing space by $r^2$. The following fluid scaled processes will be used: For $t \geq 0$,

$$\begin{aligned}
\bar{E}^r(t) &\doteq r^{-2} E^r(r^2 t), & \bar{S}^r(t) &\doteq r^{-2} S^r(r^2 t), \\
\bar{\mathbf{\Phi}}^{j,r}(t) &\doteq r^{-2} \mathbf{\Phi}^{j,r}([r^2 t]), \quad j=1,\ldots,\mathbf{J}, & \bar{T}^r(t) &\doteq r^{-2} T^r(r^2 t), \\
\bar{I}^r(t) &\doteq r^{-2} I(r^2 t), & \bar{Q}^r(t) &\doteq r^{-2} Q(r^2 t).
\end{aligned} \tag{2.19}$$

*Diffusion scaled processes.* These are obtained from the original process by accelerating time by a factor of $r^2$ and dividing space by $r$ (after appropriate centering). The following diffusion scaled processes will be considered: For $t \geq 0$,

$$\begin{aligned}
\hat{E}^r(t) &\doteq \frac{(E^r(r^2 t) - \alpha^r r^2 t)}{r}, \\
\hat{S}^r(t) &\doteq \frac{(S^r(r^2 t) - \beta^r r^2 t)}{r}, \\
\hat{\mathbf{\Phi}}^{j,r}(t) &\doteq \frac{(\mathbf{\Phi}^{j,r}([r^2 t]) - p^j [r^2 t])}{r}, \quad j = 0, 1, \ldots, \mathbf{J},
\end{aligned} \tag{2.20}$$

$$\hat{U}^r(t) \doteq r^{-1} U^r(r^2 t),$$
$$\hat{Q}^r(t) \doteq r^{-1} Q^r(r^2 t),$$
$$\hat{W}^r(t) \doteq r^{-1} W^r(r^2 t),$$
$$\hat{Y}^r(t) \doteq r^{-1} Y^r(r^2 t),$$

where for $x \in \mathbb{R}_+$, $[x]$ denotes its integer part, that is, the largest integer bounded by $x$. Note that the last four processes do not need to be centered, as it will be seen that their fluid scaled versions converge to zero as $r \to \infty$. Define, for $t \geq 0$,

$$\hat{X}_i^r(t) \doteq \hat{E}_i^r(t) - \sum_{j=1}^{\mathbf{J}} (C_{ij} - p_i^j) \hat{S}_j^r(\bar{T}_j^r(t))$$
(2.21)
$$+ \sum_{j=1}^{\mathbf{J}} \hat{\mathbf{\Phi}}_i^{j,r}(\bar{S}_j^r(\bar{T}_j^r(t))), \qquad i = 1, \ldots, \mathbf{I}.$$

Recall the definitions of $\theta_i^r$ and $\hat{q}^r$ introduced in Assumption 2.1. From the above definitions of rescaled processes, and equations (2.7), (2.8), (2.12), (2.15) and (2.17), one can establish the following relationships between the various diffusion-scaled quantities: For all $t \geq 0$,

(2.22) $\quad \hat{Q}^r(t) = \hat{q}^r + \hat{X}^r(t) + [\theta_1^r t - (C - P') \operatorname{diag}(\theta_2^r) \bar{T}^r(t)] + R \hat{Y}^r(t),$

(2.23) $\quad \hat{U}^r(t) = K \hat{Y}^r(t),$

(2.24)
$$\hat{W}^r(t) = \Lambda \hat{q}^r + \Lambda \hat{X}^r(t) + \Lambda [\theta_1^r t - (C - P') \operatorname{diag}(\theta_2^r) \bar{T}^r(t)] + \Lambda R \hat{Y}^r(t)$$
$$= \Lambda \hat{q}^r + \Lambda \hat{X}^r(t) + \Lambda [\theta_1^r t - (C - P') \operatorname{diag}(\theta_2^r) \bar{T}^r(t)] + G \hat{U}^r(t).$$

In obtaining the above equations, we have also used (2.12). Finally, in obtaining (2.24), we have used (2.18) and (2.23).

*Cost function.* For the network $\mathcal{N}^r$, we consider an expected infinite horizon discounted (linear) holding cost associated with a scheduling policy $T^r$, given in terms of the corresponding diffusion-scaled queue length process $\hat{Q}^r$ as follows:

(2.25) $\quad J^r(q^r, T^r) \doteq \mathbb{E}\left( \int_0^\infty e^{-\gamma t} h \cdot \hat{Q}^r(t) \, dt \right) + \mathbb{E}\left( \int_0^\infty e^{-\gamma t} p \cdot d\hat{U}^r(t) \right),$

where $q^r = Q^r(0)$. Here, $\gamma \in (0, \infty)$ is the "discount factor," $h$, an $\mathbf{I} \times 1$ vector with each component $h_i \in (0, \infty)$, $i = 1, \ldots, \mathbf{I}$, is the vector of "holding costs" for the $\mathbf{I}$ buffers and $p \geq 0$ is a $(\mathbf{K} + \mathbf{J} - \mathbf{B}) \times 1$ vector of costs for incurring idleness in the system.

We now introduce the key admissibility requirements on the control policy $T^r$.

*Multiparameter filtrations and stopping-times.* We refer the reader to Section 2.8 of [10] for basic definitions and properties of multiparameter filtrations, stopping times and martingales. For $m = (m_1, \ldots, m_\mathbf{I}) \in \mathbb{N}^\mathbf{I}$, $n = (n_1, \ldots, n_\mathbf{J}) \in \mathbb{N}^\mathbf{J}$ we define the multiparameter filtration generated by interarrival and service times and routing variables as

$$\begin{aligned}(2.26)\quad \mathcal{F}^r((m,n)) = \sigma\{u_i^r(m_i'), v_j^r(n_j'), \phi_i^{j,r}(n_j'): \\ m_i' \leq m_i, n_j' \leq n_j; i = 1, \ldots, \mathbf{I}, j = 1, \ldots, \mathbf{J}\}.\end{aligned}$$

Then $\{\mathcal{F}^r((m,n)) : m \in \mathbb{N}^\mathbf{I}, n \in \mathbb{N}^\mathbf{J}\}$ is a multiparameter filtration with the following (partial) ordering:

$$(2.27)\quad (m^1, n^1) \leq (m^2, n^2) \quad \text{if and only if} \quad m_i^1 \leq m_i^2, \quad n_j^1 \leq n_j^2;$$
$$i = 1, \ldots, \mathbf{I},\ j = 1, \ldots, \mathbf{J}.$$

Let

$$(2.28)\quad \mathcal{F}^r \doteq \sigma\left\{\bigcup_{(m,n) \in \mathbb{N}^{\mathbf{I}+\mathbf{J}}} \mathcal{F}^r((m,n))\right\}.$$

For all $(m, n) \in \{0, 1\}^{\mathbf{I}+\mathbf{J}}$, we define $\mathcal{F}^r((m,n)) = \mathcal{F}^r((\mathbf{1}, \mathbf{1}))$, where $\mathbf{1}$ denotes the vector of 1's.

DEFINITION 2.6 (*Admissibility of control policies*). For a fixed $r$ and $q^r \in \mathbb{R}^\mathbf{I}$, a scheduling policy $T^r = \{(T_1^r(t), \ldots, T_\mathbf{J}^r(t)) : t \geq 0\}$ is called *admissible* for $\mathcal{N}^r$ with initial condition $q^r$ if the following conditions hold:

(i) $T_j^r$ is nondecreasing, nonnegative and satisfies $T_j^r(0) = 0$ for $j = 1, \ldots, \mathbf{J}$.

(ii) $I_k^r$ defined by (2.8) is nondecreasing, nonnegative and satisfies $I_k^r(0) = 0$ for $k = 1, \ldots, \mathbf{K}$.

(iii) $Q_i^r$ defined in (2.7) is nonnegative for $i = 1, \ldots, \mathbf{I}$.

(iv) Define for each $r, t \geq 0$,

$$\begin{aligned}(2.29)\quad \sigma_0^r(t) &= (\sigma_0^{r,E}(t), \sigma_0^{r,S}(t)) \\ &\doteq (E_i^r(r^2 t) + 1 : i = 1, \ldots, \mathbf{I}; S_j^r(T_j^r(r^2 t)) + 1 : j = 1, \ldots, \mathbf{J}).\end{aligned}$$

Then, for each $t \geq 0$,

(2.30) $\quad \sigma_0^r(t)$ is a $\{\mathcal{F}^r((m,n)) : m \in \mathbb{N}^\mathbf{I}, n \in \mathbb{N}^\mathbf{J}\}$ stopping time.

Define the filtration $\{\mathcal{F}_1^r(t) : t \geq 0\}$ as

$$\mathcal{F}_1^r(t) \doteq \mathcal{F}^r(\sigma_0^r(t)) = \sigma\{A \in \mathcal{F}^r : A \cap \{\sigma_0^r(t) \leq (m,n)\} \in \mathcal{F}^r((m,n)),$$
$$\text{for all } m \in \mathbb{N}^\mathbf{I}, n \in \mathbb{N}^\mathbf{J}\},$$

for all $t \geq 0$. Then

(2.31) $\qquad\qquad\qquad \hat{U}^r$ is $\{\mathcal{F}_1^r(t)\}$ adapted.

A sequence of control policies $\{T^r\}$ is called *admissible with initial condition $q$* if, for each $r$, $T^r$ is an admissible policy for the $r$th network with some initial condition $q^r$ and $\hat{q}^r \to q$ as $r \to \infty$. Denote the class of all admissible sequences $\{T^r\}$ with initial condition $q$ by $\mathcal{A}(q)$.

REMARK 2.7. (i) and (ii) in Definition 2.6 imply, in view of (2.8) and properties of the $A$ matrix, that

$$(2.32) \quad 0 \le T_j^r(t) - T_j^r(s) \le t - s, \qquad j = 1, \ldots, \mathbf{J}, \text{ for all } 0 \le s \le t < \infty.$$

In particular, $T_j^r$ is a process with Lipschitz continuous paths. Condition (iv) in Definition 2.6 can be interpreted as a nonanticipativity condition. Indeed, in Section 5 of this paper we will show (Theorem 5.4) that a very large family of control policies that do not "anticipate the future" in a suitable sense satisfy the admissibility condition (iv) of Definition 2.6.

The following proposition gives another sufficient condition for part (iv) of the above definition to hold:

PROPOSITION 2.8. *Let $r$ and $q^r \in \mathbb{R}^{\mathbf{I}}$ be fixed and suppose that $T^r$ is a scheduling policy that, in addition to* (i)–(iii) *of Definition 2.6, satisfies, for all $t, t_j \in [0, \infty)$, $j = 1, \ldots, \mathbf{J}$ and $m \in \mathbb{N}^{\mathbf{I}}$, $n \in \mathbb{N}^{\mathbf{J}}$,*

$$(2.33) \quad \begin{aligned} \{T_j^r(t) < t_j, j = 1, \ldots, \mathbf{J}\} &\cap \{E_i^r(t) = m_i, i = 1, \ldots, \mathbf{I}\} \\ &\cap \{S_j^r(t_j) = n_j, j = 1, \ldots, \mathbf{J}\} \in \mathcal{F}^r((m, n)). \end{aligned}$$

*Then $T^r$ is admissible for $\mathcal{N}^r$ with initial condition $q^r$.*

PROOF. We first prove (2.30). For this, it suffices to show that for any $m \in \mathbb{N}^{\mathbf{I}}$, $n \in \mathbb{N}^{\mathbf{J}}$,

$$(2.34) \quad \{\sigma_{0,i}^{r,E}(t) \le m_i, \sigma_{0,j}^{r,S}(t) \le n_j, i = 1, \ldots, \mathbf{I}, j = 1, \ldots, \mathbf{J}\} \in \mathcal{F}^r((m, n)).$$

The set above can be written as

$$\left\{ \sum_{m_i'=1}^{m_i} u_i^r(m_i') > r^2 t, \sum_{n_j'=1}^{n_j} v_j^r(n_j') > T_j^r(r^2 t), i = 1, \ldots, \mathbf{I}, j = 1, \ldots, \mathbf{J} \right\}$$

$$= \bigcup_{q \in \mathbb{Q}^{\mathbf{J}}} \{T_j^r(r^2 t) < q_j, j = 1, \ldots, \mathbf{J}\}$$

$$(2.35) \qquad \cap \left\{ \sum_{m_i'=1}^{m_i} u_i^r(m_i') > r^2 t, \sum_{n_j'=1}^{n_j} v_j^r(n_j') > q_j, i = 1, \ldots, \mathbf{I}, j = 1, \ldots, \mathbf{J} \right\}$$

$$\equiv \bigcup_{q \in \mathbb{Q}^{\mathbf{J}}} \mathcal{A}_q \cap \mathcal{B}_q$$

$$= \bigcup_{q \in \mathbb{Q}^{\mathbf{J}}} \bigcup_{(\tilde{m},\tilde{n}) < (m,n)} \mathcal{A}_q \cap \mathcal{B}_q \cap \{E^r(r^2 t) = \tilde{m}, S^r(q) = \tilde{n}\},$$

where $\mathbb{Q}$ denotes the set of rational numbers. From (2.33), it follows that for each $q$ and $(\tilde{m},\tilde{n}) < (m,n)$ the set $\mathcal{A}_q \cap \{E^r(r^2 t) = \tilde{m}, S^r(q) = \tilde{n}\} \in \mathcal{F}^r((\tilde{m},\tilde{n})) \subset \mathcal{F}^r((m,n))$. Also, clearly, $\mathcal{B}_q \in \mathcal{F}^r((m,n))$ for every $q \in \mathbb{Q}^{\mathbf{J}}$. Thus, we have that the set in (2.35) is in $\mathcal{F}^r((m,n))$ and hence (2.34) follows.

Since $\sigma_0^r(t)$ is a nondecreasing function of $t$, it follows that $\{\mathcal{F}_1^r(t) : t \geq 0\}$ is a filtration. In order to prove (2.31), it suffices to show that $T_j^r(r^2 \cdot)$ is adapted to the filtration $\{\mathcal{F}_1^r(t) : t \geq 0\}$. For this, it is enough to show that for fixed $t \geq 0$ and real $a_j$, $j = 1, \ldots, \mathbf{J}$, we have $\{T_j^r(r^2 t) < a_j, j = 1, \ldots, \mathbf{J}\} \in \mathcal{F}_1^r(t)$. The last statement will follow if the following property holds for all $m, n$:

$$\{T_j^r(r^2 t) < a_j, j = 1, \ldots, \mathbf{J}\}$$
$$\cap \{\sigma_{0,i}^{r,E}(t) \leq m_i, \sigma_{0,j}^{r,S}(t) \leq n_j, i = 1, \ldots, \mathbf{I}, j = 1, \ldots, \mathbf{J}\} \in \mathcal{F}^r((m,n)).$$

The set above can be written as

$$\bigcup_{q \in \mathbb{Q}^{\mathbf{J}}} \{T_j^r(r^2 t) < q_j, j = 1, \ldots, \mathbf{J}\}$$

$$\cap \left\{ \sum_{m_i'=1}^{m_i} u_i^r(m_i') > r^2 t, \left( \sum_{n_j'=1}^{n_j} v_j^r(n_j') \wedge a_j \right) > q_j, i = 1, \ldots, \mathbf{I}, j = 1, \ldots, \mathbf{J} \right\}.$$

Finally, proceeding exactly as in the proof of the first part, one obtains that this set is in $\mathcal{F}^r((m,n))$. This completes the proof of the proposition. □

Given a sequence of admissible control policies $\{T^r\}$, with initial condition sequence $\{q^r\}$ such that $\hat{q}^r \to q$ as $r \to \infty$, the associated asymptotic cost is defined as

(2.36) $$J(q, \{T^r\}) \doteq \liminf_{r \to \infty} J^r(q^r, T^r).$$

Let the infimum of the asymptotic cost over all admissible sequences $\{T^r\}$ be denoted by $J^*(q)$, that is,

(2.37) $$J^*(q) \doteq \inf_{\{T^r\} \in \mathcal{A}(q)} J(q, \{T^r\}).$$

The main goal of the paper is to give a lower bound for this infimum in terms of the value function of an associated diffusion control problem.

*Associated diffusion control problems.* We now introduce the formal diffusion approximation to the network control problem described above. It will be seen in Lemma 3.3 [part (C)], via an application of a functional central limit theorem, that $\hat{X}^r$ defined in (2.21) converges weakly to a Brownian motion $\tilde{X}$ with drift 0 and covariance matrix

$$(2.38) \quad \Sigma \doteq \Sigma^u + (C - P')\Sigma^v \operatorname{diag}(x^*)(C - P')' + \sum_{j=1}^{\mathbf{J}}[\Sigma^{\phi^j}\beta_j x_j^*],$$

where $\Sigma^u$ is an $\mathbf{I} \times \mathbf{I}$ diagonal matrix with diagonal entries $(\sigma_i^u)^2$, $i = 1,\ldots,\mathbf{I}$, $\Sigma^v$ is a $\mathbf{J} \times \mathbf{J}$ diagonal matrix with diagonal entries $(\sigma_j^v)^2$, $j = 1,\ldots,\mathbf{J}$ and the $\Sigma^{\phi^j}$'s are $\mathbf{I} \times \mathbf{I}$ matrices with entries $\sigma^{\phi^j}_{i_1 i_2}, i_1, i_2 = 1,\ldots,\mathbf{I}$ [cf. (2.3)]. Furthermore, it will be seen in Lemma 3.3 [part (B)] that for all "reasonable" control policies, the third term on the right-hand side of (2.22) converges to $\theta t$, where

$$(2.39) \quad \theta \doteq \theta_1 - (C - P')\operatorname{diag}(\theta_2)x^*.$$

Although the weak convergence of the last term (indeed, even its tightness) in the right-hand side of (2.22) is far from obvious on formal passing to the limit, one is led to the following diffusion control problem:

DEFINITION 2.9 (*Brownian control problem* (*BCP*)). A **J**-dimensional adapted process $\tilde{Y}$, defined on some filtered probability space $(\tilde{\Omega}, \tilde{\mathcal{F}}, \tilde{\mathbb{P}}, \{\tilde{\mathcal{F}}(t)\})$ which supports an **I**-dimensional $\{\mathcal{F}(t)\}$-Brownian motion $\tilde{X}$ with drift 0 and covariance matrix $\Sigma$ defined in (2.38), is called an *admissible control* for the Brownian control problem with the initial condition $q$ iff the following two properties hold $\tilde{\mathbb{P}}$-a.s.:

$$\tilde{Q}(t) \doteq q + \theta t + \tilde{X}(t) + R\tilde{Y}(t) \geq 0 \quad \text{for all } t \geq 0,$$
$$\tilde{U}(t) \doteq K\tilde{Y}(t) \text{ is nondecreasing and } \tilde{U}(0) \geq 0,$$

where $\theta$ is as in (2.39). We refer to $\Gamma = (\tilde{\Omega}, \tilde{\mathcal{F}}, \tilde{\mathbb{P}}, \{\tilde{\mathcal{F}}(t)\}, \tilde{X})$ as a *system* and $\tilde{X}$ as a $(0, \Sigma)$-*Brownian motion*. We denote the class of all such admissible controls by $\tilde{\mathcal{A}}(q)$. The Brownian control problem is to infimize

$$(2.40) \quad \tilde{J}(\tilde{Y}) \doteq \tilde{\mathbb{E}}\int_0^\infty e^{-\gamma t} h \cdot \tilde{Q}(t)\,dt + \tilde{\mathbb{E}}\int_{[0,\infty)} e^{-\gamma t} p \cdot d\tilde{U}(t)$$

over all admissible controls $\tilde{Y} \in \tilde{\mathcal{A}}(q)$. Define the value function

$$(2.41) \quad \tilde{J}^*(q) = \inf_{\tilde{Y} \in \tilde{\mathcal{A}}(q)} \tilde{J}(\tilde{Y}).$$

Since an admissible control is not required to be of bounded variation, the above control problem is somewhat nonstandard and is difficult to analyze directly. However, as shown in [11], under suitable conditions, one can replace this control problem by an equivalent, more classical control problem which is also typically of a much lower state dimension. This control problem, in stochastic control terminology, is a problem of *singular control with state constraints* (SCSC) and is given as follows. [In the literature, this problem is also referred to as the *equivalent workload formulation* (EWF) of the BCP.]

*Effective cost function.* The definition of the equivalent workload formulation (EWF) is tied to the notion of an effective cost function, described as follows. Recall the definition of the workload matrix $\Lambda$ introduced in (2.17). Let $\mathcal{W} \doteq \{\Lambda z : z \in \mathbb{R}_+^{\mathbf{I}}\}$. For each $w \in \mathcal{W}$, define

$$\hat{h}(w) \doteq \min\{h \cdot q : \Lambda q = w, q \geq 0\}. \tag{2.42}$$

It is well known (see Theorem 2 of [4]) that one can take a continuous selection of the minimizer in the above linear program. That is, there is a continuous map $q^* : \mathcal{W} \to \mathbb{R}_+^{\mathbf{I}}$ such that

$$q^*(w) \in \arg\min\{h \cdot q : \Lambda q = w, q \geq 0\}. \tag{2.43}$$

As an immediate consequence of these definitions, if $w(q) \doteq \Lambda q, q \in \mathbb{R}_+^{\mathbf{I}}$ and $w_r, w \in \mathcal{W}$, then

$$\begin{aligned} h \cdot q &\geq \hat{h}(w(q)) \qquad \text{for all } q \in \mathbb{R}_+^{\mathbf{I}}, \\ \hat{h}(w_r) &\to \hat{h}(w) \qquad \text{whenever } w_r \to w \text{ as } r \to \infty. \end{aligned} \tag{2.44}$$

The function $\hat{h}(\cdot)$ will be referred to as the *effective cost function*. The equivalent workload formulation (EWF) and the associated control problem is defined below. Let

$$\mathcal{K} \doteq \{u \in \mathbb{R}^{\mathbf{K}+\mathbf{J}-\mathbf{B}} | u = Ky, y \in \mathbb{R}^{\mathbf{J}}\}. \tag{2.45}$$

DEFINITION 2.10 (*Equivalent workload formulation* (*EWF*)). A ($\mathbf{K} + \mathbf{J} - \mathbf{B}$)-dimensional adapted process $\tilde{U}_0$, defined on some filtered probability space $(\tilde{\Omega}, \tilde{\mathcal{F}}, \tilde{\mathbb{P}}, \{\tilde{\mathcal{F}}(t)\})$ which supports an $\mathbf{I}$-dimensional $\{\mathcal{F}(t)\}$-Brownian motion $\tilde{X}$ with drift 0 and covariance matrix $\Sigma$ defined in (2.38), is called an *admissible control* for the equivalent workload formulation of the Brownian control problem with the initial condition $q$ iff the following two properties hold $\tilde{\mathbb{P}}$-a.s.:

$$\begin{aligned} &\tilde{U}_0 \text{ is nondecreasing, } \tilde{U}_0(0) \geq 0, \tilde{U}_0(t) \in \mathcal{K} \qquad \text{for all } t \geq 0, \\ &\tilde{W}(t) \doteq w + \Lambda \theta t + \Lambda \tilde{X}(t) + G\tilde{U}_0(t) \in \mathcal{W} \qquad \text{for all } t \geq 0, \end{aligned} \tag{2.46}$$

where $\theta$ is as in (2.39). We denote the class of all such admissible controls by $\tilde{\mathcal{A}}_0(w)$. The control problem for the equivalent workload formulation is to infimize

$$(2.47) \qquad \tilde{J}_0(\tilde{U}_0) \doteq \tilde{\mathbb{E}} \int_0^\infty e^{-\gamma t} \hat{h}(\tilde{W}(t))\, dt + \tilde{\mathbb{E}} \int_{[0,\infty)} e^{-\gamma t} p \cdot d\tilde{U}_0(t)$$

over all admissible controls $\tilde{U}_0 \in \tilde{\mathcal{A}}_0(w)$. Define the value function

$$(2.48) \qquad \tilde{J}_0^*(w) = \inf_{\tilde{U}_0 \in \tilde{\mathcal{A}}_0(w)} \tilde{J}_0(\tilde{U}_0).$$

From Theorem 2 of [11], it follows that for all $w \in \mathcal{W}, q \in \mathbb{R}_+^{\mathbf{I}}$ satisfying $w = \Lambda q$,

$$(2.49) \qquad \tilde{J}^*(q) = \tilde{J}_0^*(w).$$

**3. Main result and proofs.** In this section, we present the main theorem of this paper (Theorem 3.1) and provide proofs of the key results that are used in proving the theorem.

THEOREM 3.1. *Let $\{T^r\}$ be an admissible sequence of control policies for the initial condition $q$, that is, for each $r$, $T^r$ is an admissible policy for $\mathcal{N}^r$ with some initial condition $q^r$, and $\hat{q}^r \to q$ as $r \to \infty$. Let $J(q, \{T^r\}) \doteq \liminf_{r \to \infty} J^r(q^r, T^r)$ be the associated asymptotic cost. Let $w \doteq \Lambda q$ and $\tilde{J}_0^*(w)$ be, as in (2.48), the value function in the EWF. Then*

$$(3.1) \qquad J(q, \{T^r\}) \geq \tilde{J}_0^*(w).$$

As a corollary of Theorem 3.1 and Theorem 2 of [11] [see (2.49)], the infimum of the cost in the network control problem is asymptotically bounded below by the value function of the BCP:

COROLLARY 3.2. *Fix $q \in \mathbb{R}_+^{\mathbf{I}}$ and set $w \doteq \Lambda q$. Let $J^*(q)$, $\tilde{J}^*(q)$ and $\tilde{J}_0^*(w)$ be given via (2.37), (2.41) and (2.48), respectively. Then*

$$(3.2) \qquad J^*(q) \geq \tilde{J}_0^*(w) = \tilde{J}^*(q).$$

The rest of the paper is devoted to the proof of Theorem 3.1. Thus, henceforth we will fix an admissible sequence $\{T^r\} \in \mathcal{A}(q)$ such that $T^r$ is an admissible policy for the initial condition $q^r$ and $\hat{q}^r \to q$ as $r \to \infty$. Furthermore, we will assume without loss of generality that

$$(3.3) \qquad J(q, \{T^r\}) \doteq \liminf_{r \to \infty} J^r(q^r, T^r) < \infty,$$

since otherwise inequality (3.1) holds trivially. Finally, the subsequence $\{r'\} \subset \{r\}$ along which the above lim inf is attained will be relabeled once more as $\{r\}$. In particular, by this relabeling,

$$J(q, \{T^r\}) = \lim_{r \to \infty} J^r(q^r, T^r). \tag{3.4}$$

We begin with some preliminary results.

LEMMA 3.3. *(A) The following functional central limit theorem holds:*

$$(\hat{E}^r, \hat{S}^r, \hat{\boldsymbol{\Phi}}^{r,j}, j = 1, \ldots, \mathbf{J}) \tag{3.5}$$
$$\Rightarrow (\tilde{E}, \tilde{S}, \tilde{\boldsymbol{\Phi}}^j, j = 1, \ldots, \mathbf{J}) \qquad \text{as } r \to \infty,$$

*where $\tilde{E}, \tilde{S}, \tilde{\boldsymbol{\Phi}}^j$, $j = 1, \ldots, \mathbf{J}$ are independent Brownian motions with drift zero and covariances $\Sigma^u, \Sigma^v, \Sigma^{\phi_j}$, $j = 1, \ldots, \mathbf{J}$ [cf. (2.38)], respectively.*

*(B) The following functional law of large numbers holds:*

$$(\bar{Q}^r, \bar{I}^r, \bar{T}^r, \bar{E}^r, \bar{S}^r, \bar{\boldsymbol{\Phi}}^{j,r}, \; j = 1, \ldots, \mathbf{J}) \tag{3.6}$$
$$\Rightarrow (\mathbf{0}, \mathbf{0}, \bar{T}^*, \alpha e, \beta e, p^j e, \; j = 1, \ldots, \mathbf{J}) \qquad \text{as } r \to \infty,$$

*where $\mathbf{0}$ denotes the process that equals 0 a.s. at all times, $e(t) \doteq t$, $\bar{T}^*(t) \doteq x^* t, t \geq 0$ and $x^*$ is as in Remark 2.4.*

*(C) Let $\hat{X}^r$ be as in (2.21) and define*

$$\tilde{X}_j(\cdot) \doteq \tilde{E}_j(\cdot) - \sum_j (C_{ij} - P'_{ij}) \tilde{S}_j(x_j^* \cdot) + \sum_j \hat{\boldsymbol{\Phi}}_i^j(x_j^* \beta_j \cdot). \tag{3.7}$$

*Then $\tilde{X}$ is a $(0, \Sigma)$-Brownian motion and*

$$\hat{\Xi}^r \doteq (\hat{E}_i^r(\cdot), \hat{S}_j^r(\bar{T}_j^r(\cdot)), \hat{\boldsymbol{\Phi}}_i^{j,r}(\bar{S}_j^r(\bar{T}_j^r(\cdot))), \hat{X}_i^r(\cdot) : i = 1, \ldots, \mathbf{I}, j = 1, \ldots, \mathbf{J}) \tag{3.8}$$
$$\Rightarrow (\tilde{E}_i(\cdot), \tilde{S}_j(x_j^* \cdot), \tilde{\boldsymbol{\Phi}}_i^j(x_j^* \beta_j \cdot) \tilde{X}_i(\cdot) : i = 1, \ldots, \mathbf{I}, j = 1, \ldots, \mathbf{J}) \doteq \Xi.$$

PROOF. Weak convergence of $\hat{E}^r, \hat{S}^r$ to Brownian motions follows from the functional central limit theorem for renewal processes (see Theorem 14.6 of [3]), and the weak convergence of $\hat{\boldsymbol{\Phi}}^{j,r}$ to a Brownian motion with the stated covariance matrix follows from Donsker's theorem (Theorem 8.2 of [3]). Finally, the joint weak convergence in (3.5) and claimed mutual independence follow from the mutual independence of the primitive processes on the left-hand side of (3.5). This proves part (A).

Part (C) is an immediate consequence of (A) and (B) on using Lemma 3.14.1 of [3]. Finally, for the proof of (B), note that from (3.5) it follows that

$$(\bar{E}^r, \bar{S}^r, \bar{\boldsymbol{\Phi}}^{j,r}, \; j = 1, \ldots, \mathbf{J}) \tag{3.9}$$
$$\Rightarrow (\alpha e, \beta e, p^j e, \; j = 1, \ldots, \mathbf{J}) \qquad \text{as } r \to \infty.$$

Also, from Definition 2.6, we have that for each $r$, $T^r$ is Lipschitz continuous with Lipschitz constant less than or equal to 1, and that $T^r(t) \leq t$ for all $t \geq 0$ (see Remark 2.7). Thus the same properties hold for the fluid-scaled processes $\bar{T}^r$ for each $r$. This shows that the sequence $\{\bar{T}^r\}$ is $\mathcal{C}$-tight. In addition, we have from (2.7), (2.8) and the definition of fluid-scaled processes in (2.19) that

$$(3.10) \quad \bar{Q}_i^r(t) = \frac{q^r}{r^2} + \bar{E}_i^r(t) - \sum_{j=1}^{\mathbf{J}} C_{ij} \bar{S}_j^r(\bar{T}_j^r(t)) + \sum_{j=1}^{\mathbf{J}} \bar{\Phi}_i^{j,r}(\bar{S}_j^r(\bar{T}_j^r(t))),$$

$$(3.11) \quad \bar{I}_k^r(t) = t - \sum_{j=1}^{\mathbf{J}} A_{kj} \bar{T}_j^r(t),$$

for all $i = 1, \ldots, \mathbf{I}$, $k = 1, \ldots, \mathbf{K}$. Combining the above representations with (3.9) and the $\mathcal{C}$-tightness of $\{\bar{T}^r\}$, we see that

$$(3.12) \quad \{(\bar{Q}^r, \bar{I}^r, \bar{T}^r, \bar{E}^r, \bar{S}^r, \bar{\Phi}^{j,r}, \ j = 1, \ldots, \mathbf{J})\} \text{ is } \mathcal{C}\text{-tight.}$$

Hence, it is enough to prove that any weak limit point of the sequence in (3.12) is given by the right-hand side of (3.6). Consider a further subsequence (we use the same subscript $r$ to denote this subsubsequence) which converges weakly to some $(\bar{Q}, \bar{I}, \bar{T}, \bar{E}, \bar{S}, \bar{\Phi}^j, \ j = 0, 1, \ldots, \mathbf{J})$. Using the Skorohod representation theorem, we can assume without loss of generality that

$$(3.13) \quad \begin{aligned} &(\bar{Q}^r, \bar{I}^r, \bar{T}^r, \bar{E}^r, \bar{S}^r, \bar{\Phi}^{j,r}, \ j = 1, \ldots, \mathbf{J}) \\ &\to (\bar{Q}, \bar{I}, \bar{T}, \bar{E}, \bar{S}, \bar{\Phi}^j, \ j = 1, \ldots, \mathbf{J}) \qquad \text{u.o.c. as } r \to \infty, \end{aligned}$$

with probability 1. Since we have $J(q, \{T^r\}) < \infty$, using Fatou's lemma we get

$$(3.14) \quad \begin{aligned} 0 &= \lim_{r \to \infty} \frac{J^r(q^r, T^r)}{r} \\ &\geq \liminf_{r \to \infty} \mathbb{E}\left[\int_0^\infty e^{-\gamma t} h \cdot \bar{Q}^r(t) \, dt\right] \\ &\geq \mathbb{E}\left[\int_0^\infty e^{-\gamma t} \left(\liminf_{r \to \infty} h \cdot \bar{Q}^r(t)\right) dt\right] \\ &= \mathbb{E}\left[\int_0^\infty e^{-\gamma t} h \cdot \bar{Q}(t) \, dt\right]. \end{aligned}$$

Recalling that $h > 0$, we have, in view of the path continuity of $\bar{Q}$, that $\bar{Q} = \mathbf{0}$, a.s. Hence, from (3.9), (3.10), (3.11), (3.12) and the definition of $R$ in (2.10), we obtain that

$$(3.15) \quad \mathbf{0} = \alpha e - R\bar{T}, \qquad \mathbf{0} \leq \bar{I} = e - A\bar{T}.$$

Finally, combining (3.15) with Assumption 2.3 [see Definition 2.2(ii) and Remark 2.4] yields

$$\bar{T}(t) = x^* t, \qquad \bar{I}(t) = et - Ax^* t = \mathbf{0}, \qquad t \geq 0. \tag{3.16}$$

This completes the proof of part (B). □

The following corollary is a direct consequence of Lemma 3.3:

COROLLARY 3.4. *Let*

$$\hat{M}^r(\cdot) \doteq (\hat{M}_i^{E,r}(\cdot) : i = 1, \ldots, \mathbf{I}, \hat{M}_j^{S,r}(\cdot) : j = 1, \ldots, \mathbf{J}, \tag{3.17}$$

$$\hat{M}_i^{\Phi,r}(\cdot) : i = 1, \ldots, \mathbf{I}),$$

*where for $i = 1, \ldots, \mathbf{I}$, $j = 1, \ldots, \mathbf{J}$ and $t \geq 0$,*

$$\hat{M}_i^{E,r}(t) \doteq \hat{E}_i^r(t), \qquad \hat{M}_j^{S,r}(t) = \hat{S}_j^r(\bar{T}_j^r(t)), \tag{3.18}$$

$$\hat{M}_i^{\Phi,r}(t) = \sum_{j=1}^{\mathbf{J}} \hat{\Phi}_i^{j,r}(\bar{S}_j^r(\bar{T}_j^r(t))).$$

*Then*

$$\hat{M}^r(\cdot) \Rightarrow M'(\cdot) \doteq \left( \tilde{E}(\cdot), i = 1, \ldots, \mathbf{I}, \tilde{S}_j(x_j^* \cdot), \ j = 1, \ldots, \mathbf{J}, \right.$$

$$\left. \sum_{j=1}^{\mathbf{J}} \tilde{\Phi}_i^j(\beta_j x_j^* \cdot), \ i = 1, \ldots, \mathbf{I} \right), \tag{3.19}$$

*where $\tilde{E}, \tilde{S}$ and $\tilde{\Phi}$ are as in Lemma 3.3(A).*

LEMMA 3.5. *There exists a $c^* \in (0, \infty)$ such that for all $i = 1, \ldots, \mathbf{I}$, $j = 1, \ldots, \mathbf{J}$, $r \geq 1$ and $t \geq 0$,*

$$\mathbb{E}\left[ \sup_{0 \leq s \leq t} |\hat{E}_i^r(s)|^2 \right] < c^*(t+1),$$

$$\mathbb{E}\left[ \sup_{0 \leq s \leq t} |\hat{S}_j^r(\bar{T}_j^r(s))|^2 \right] < c^*(t+1), \tag{3.20}$$

$$\mathbb{E}\left[ \sup_{0 \leq s \leq t} |\hat{\Phi}_i^{j,r}(\bar{S}_j^r(\bar{T}_j^r(s)))|^2 \right] < c^*(t+1).$$

*Furthermore, for $m = 1, \ldots, \mathbf{K} + \mathbf{J} - \mathbf{B}$, $t \geq 0$,*

$$\limsup_{r \to \infty} \mathbb{E}[\hat{U}_m^r(t)] < \infty. \tag{3.21}$$

PROOF. We first consider the third bound in (3.20). Fix $i \in \{1, \ldots, \mathbf{I}\}$, $j \in \{1, \ldots, \mathbf{J}\}$. Note that

$$
\begin{aligned}
(3.22) \quad \sup_{0 \leq s \leq t} |\hat{\mathbf{\Phi}}_i^{j,r}(\bar{S}_j^r(\bar{T}_j^r(s)))| &\leq \sup_{0 \leq s \leq t} |\hat{\mathbf{\Phi}}_i^{j,r}(\bar{S}_j^r(s))| \\
&\leq \frac{2}{r} + \frac{1}{r}\left[\sup_{0 \leq s \leq t} |\tilde{\zeta}_i^{j,r}(\sigma_j^{r,S}(s))|\right],
\end{aligned}
$$

where $\tilde{\zeta}_i^{j,r}(n) \doteq \sum_{k=1}^n (\phi_i^{j,r}(k) - p_i^j)$, $n = 1, 2, \ldots$ and $\sigma_j^{r,S}(s) = S_j^r(r^2 s) + 1$, $s \geq 0$. The above bound follows from the fact that $\bar{T}_j^r(s) \leq s$, $s \geq 0$, and that $|\phi_i^{j,r}(k) - p_i^j| \leq 2$, for all $k = 1, 2, \ldots$. Define

$$\tilde{\mathcal{F}}(n) \doteq \{\phi_i^{j,r}(k), v_j^r(k) : k \leq n\}, \qquad n \geq 1.$$

It is easy to check that $\{\tilde{\zeta}_i^{j,r}(n)\}$ is a martingale and $\sigma_j^{r,S}(s)$, for each $s \geq 0$, is a stopping time with respect to the filtration $\{\tilde{\mathcal{F}}(n)\}$. Furthermore, the quadratic variation of this martingale is given as

$$(3.23) \qquad \langle \tilde{\zeta}_i^{j,r}\rangle(n) \doteq (\sigma_{ii}^{\phi^j})n, \qquad n = 1, 2, \ldots.$$

Next, note that for all $n = 1, 2, \ldots$,

$$(3.24) \quad \mathbb{E}[|\tilde{\zeta}_i^{j,r}(n+1) - \tilde{\zeta}_i^{j,r}(n)||\tilde{\mathcal{F}}(n)] = \mathbb{E}|\phi_i^{j,r}(n+1) - p_i^j| \leq 2,$$

and that $\{\sigma_j^{r,S}(s) : s \geq 0\}$ is a nondecreasing sequence of stopping times, satisfying

$$
\begin{aligned}
(3.25) \quad \mathbb{E}[\sigma_j^{r,S}(s)] &= \mathbb{E}(S_j^r(r^2 s) + 1) \\
&\leq c_1(r^2 s + 1) < \infty, \qquad s \geq 0,
\end{aligned}
$$

where $c_1 \in (0, \infty)$ is a constant independent of $r$ and $s$. This inequality is a consequence of Assumption 2.1 and well-known moment inequalities for renewal processes. Thus, using an optional sampling theorem (see Theorems 4.7.4, 4.7.5 of [9]) we obtain that $\{\tilde{\zeta}_i^{j,r}(\sigma_j^{r,S}(s)), \tilde{\mathcal{F}}(\sigma_j^{r,S}(s)) : s \geq 0\}$ is a martingale. Using (3.23) and (3.25), we have, for all $t \geq 0$,

$$
\begin{aligned}
(3.26) \quad \mathbb{E}[|\tilde{\zeta}_i^{j,r}(\sigma_j^{r,S}(t))|^2] &= \sigma_{ii}^{\phi^j} \mathbb{E}[\sigma_j^{r,S}(t)] \\
&\leq c_1 \sigma_{ii}^{\phi^j}(r^2 t + 1).
\end{aligned}
$$

Taking the expectation of the squares of both sides of (3.22) and using Doob's inequality, we obtain that

$$\mathbb{E}\left[\sup_{0 \leq s \leq t} |\hat{\mathbf{\Phi}}_i^{j,r}(\bar{S}_j^r(\bar{T}_j^r(s)))|^2\right]$$

$$
\begin{align}
(3.27) \qquad &\leq 2\left[\frac{4}{r^2} + \frac{1}{r^2}\mathbb{E}\left(\sup_{0\leq s\leq t}|\tilde{\zeta}_i^{j,r}(\sigma_j^{r,S}(s))|^2\right)\right] \\
&\leq \frac{8}{r^2}(1 + \mathbb{E}|\tilde{\zeta}_i^{j,r}(\sigma_j^{r,S}(t))|^2).
\end{align}
$$

Combining the above display with (3.26), we obtain the third bound in (3.20).

We will next show that for some $c^* \in (0,\infty)$, we have for all $j=1,\ldots,\mathbf{J}$, $r>1$ and $t\geq 0$ that

$$
(3.28) \qquad \mathbb{E}\left[\sup_{0\leq s\leq t}|\hat{S}_j^r(s)|^2\right] < c^*(t+1).
$$

Proof for the first term in (3.20) is almost identical to the proof of (3.28), and the second bound in (3.20) follows from (3.28) using the fact that $\bar{T}_j^r(t) \leq t$ for all $t\geq 0$. Define for $t\geq 0$, $n = 1,2,\ldots$

$$
(3.29) \quad \tilde{\eta}_j^r(n) \doteq \frac{1}{r}\sum_{k=1}^n (1 - \beta_j^r v_j^r(k)), \qquad \varepsilon_j^{S,r}(t) \doteq \sup_{0\leq s\leq t}|\hat{S}_j^r(s) - \tilde{\eta}_j^r(\sigma_j^{r,S}(s))|.
$$

As a convention, we set $\tilde{\eta}_j^r(0) \equiv 0$. A straightforward calculation shows that

$$
\begin{align}
\varepsilon_j^{S,r}(t) &\leq \frac{\beta_j^r}{r}\left[\sup_{0\leq s\leq t} v_j^r(\sigma_j^{r,S}(s))\right] + \frac{1}{r} \\
(3.30) \qquad &\leq \left[\sup_{0\leq s\leq t}|\tilde{\eta}_j^r(S_j^r(r^2s)+1) - \tilde{\eta}_j^r(S_j^r(r^2s))|\right] + \frac{2}{r} \\
&\leq 2\left[\sup_{0\leq s\leq t}|\tilde{\eta}_j^r(S_j^r(r^2s)+1)|\right] + \frac{2}{r}.
\end{align}
$$

Hence, from (3.29) we have that

$$
\begin{align}
\sup_{0\leq s\leq t}|\hat{S}_j^r(s)| &\leq \sup_{0\leq s\leq t}|\tilde{\eta}_j^r(\sigma_j^{r,S}(s))| + \varepsilon_j^{S,r}(t) \\
(3.31) \qquad &\leq 3\left[\sup_{0\leq s\leq t}|\tilde{\eta}_j^r(\sigma_j^{r,S}(s))|\right] + \frac{2}{r}.
\end{align}
$$

Note that $\{\tilde{\eta}_j^r(n)\}$ is a square integrable $\{\tilde{\mathcal{F}}(n)\}$-martingale and for that each $s\geq 0$, $\sigma_j^{r,S}(s)$ is a $\{\tilde{\mathcal{F}}(n)\}$-stopping time. Hence, an argument analogous to the one leading to (3.27) shows that

$$
\begin{align}
\mathbb{E}\left[\sup_{0\leq s\leq t}|\tilde{\eta}_j^r(\sigma_j^{r,S}(s))|^2\right] &\leq 4\mathbb{E}[|\tilde{\eta}_j^r(\sigma_j^{r,S}(s))|^2] \\
(3.32) \qquad &\leq 4\frac{(\beta_j^r \sigma_j^{v,r})^2}{r^2}\mathbb{E}[S_j^r(r^2t)+1] \\
&\leq \frac{c_2}{r^2}(r^2t+1),
\end{align}
$$

where, once again, $c_2$ is a constant independent of $r$ and $t$. Combining this estimate with (3.31) proves (3.20).

Finally, we consider (3.21). Since $h$ is strictly positive, one can find $c_3 \in (0, \infty)$ such that

$$|\hat{W}^r(t)| = |\Lambda \hat{Q}^r(t)| \leq |\Lambda| |\hat{Q}^r(t)| \leq c_3 h \cdot \hat{Q}^r(t), \qquad r \geq 1 \text{ and } t \geq 0.$$

Thus, from (3.3) we get

$$(3.33) \qquad \limsup_{r \to \infty} \mathbb{E} \int_0^\infty e^{-\gamma s} |\hat{W}^r(s)| \, ds < \infty.$$

Furthermore, from (2.24), (3.20) and Assumption 2.5, one can find $c_4 \in (0, \infty)$ such that

$$c \mathbb{E} \hat{U}_m^r(t) \leq c \mathbb{E} |\hat{U}^r(t)| \leq \mathbb{E} |G \hat{U}^r(t)| \leq \mathbb{E} |\hat{W}^r(t)| + c_4(t+1),$$

for all $m = 1, \ldots, \mathbf{K} + \mathbf{J} - \mathbf{B}$, $r \geq 1$ and $t \geq 0$. Combining this with (3.33), we have

$$(3.34) \qquad \limsup_{r \to \infty} \mathbb{E} \int_0^\infty e^{-\gamma s} \hat{U}_m^r(s) \, ds < \infty.$$

Finally, using the monotonicity of $\hat{U}_m^r$, we get, for all $t \geq 0$,

$$\mathbb{E} \int_0^\infty e^{-\gamma s} \hat{U}_m^r(s) \, ds \geq \mathbb{E} \int_t^{t+1} e^{-\gamma s} \hat{U}_m^r(s) \, ds \geq e^{-\gamma(t+1)} \mathbb{E}(\hat{U}_m^r(t)).$$

Combining this inequality with (3.34), we have (3.21). $\square$

In the following lemma, we introduce the time transformation which plays a crucial role in the proof of Theorem 3.1:

LEMMA 3.6. *For each $r$, define the process $\tau^r(\cdot)$ with values in $[0, \infty)$ as follows:*

$$(3.35) \qquad \tau^r(t) \doteq t + \sum_{m=1}^{\mathbf{K}+\mathbf{J}-\mathbf{B}} \hat{U}_m^r(t), \qquad t \geq 0.$$

*Then the map $\tau^r : [0, \infty) \to [0, \infty)$ is (almost surely) continuous and strictly increasing. Define $\check{\tau}^r(t) \doteq \inf\{s \geq 0 : \tau^r(s) > t\}$, $t \geq 0$. Consider the following time-transformed processes:*

$$(3.36) \qquad \begin{aligned} \check{W}^r(\cdot) &\doteq \hat{W}^r(\check{\tau}^r(\cdot)), & \check{X}^r(\cdot) &\doteq \hat{X}^r(\check{\tau}^r(\cdot)), \\ \check{U}^r(\cdot) &\doteq \hat{U}^r(\check{\tau}^r(\cdot)), & \check{M}^r(\cdot) &\doteq \hat{M}^r(\check{\tau}^r(\cdot)). \end{aligned}$$

*Then*

(A) $\{(\check{W}^r, \check{X}^r, \check{U}^r, \check{\tau}^r, \hat{\Xi}^r, \hat{M}^r, \check{M}^r)\}_r$ *is $\mathcal{C}$-tight.*

(B) *Let a weak limit point of the above sequence, $(\check{W}(\cdot),\check{X}(\cdot),\check{U}(\cdot),\check{\tau}(\cdot),\Xi(\cdot),$ $M'(\cdot),\check{M}(\cdot))$, be given on some probability space $(\hat{\Omega},\hat{\mathcal{F}},\hat{\mathbb{P}})$. Then, almost surely,*

$$\check{\tau}(t) \to \infty \quad \text{as } t \to \infty \tag{3.37}$$

*Define $\tau(t) = \inf\{s \geq 0: \check{\tau}(s) > t\}$, $t \geq 0$. With probability one, $\check{\tau}$ is continuous, $\tau$ is right continuous and both $\tau,\check{\tau}$ are nondecreasing maps from $[0,\infty)$ to itself.*

(C) *The process $M'$ equals the right-hand side of* (3.19), *with $\tilde{E},\tilde{S}$ and $\tilde{\boldsymbol{\Phi}}$ as in Lemma* 3.3(A). *Furthermore, the following equality holds almost surely:*

$$\begin{aligned}
\check{M}(\cdot) &\equiv (\check{M}^E(\cdot), \check{M}^S(\cdot), \check{M}^\Phi(\cdot)) \\
&= \bigg( \tilde{E}(\check{\tau}(\cdot)),\ i=1,\ldots,\mathbf{I},\ \tilde{S}_j(x_j^*\check{\tau}(\cdot)),\ j=1,\ldots,\mathbf{J}, \\
&\qquad \sum_{j=1}^{\mathbf{J}} \tilde{\boldsymbol{\Phi}}_i^j(\beta_j x_j^* \check{\tau}(\cdot)),\ i=1,\ldots,\mathbf{I} \bigg).
\end{aligned} \tag{3.38}$$

PROOF. From Remark 2.7, we have that $0 \leq T_j^r(t) - T_j^r(s) \leq t-s$ for all $0 \leq s \leq t$ and $j=1,\ldots,\mathbf{J}$. This, along with (2.14) and (2.23), shows that $\tau^r$ has Lipschitz continuous paths. Also, since $\hat{U}_m^r$ is nonnegative, $\tau^r$ is strictly increasing. This proves the first statement in the lemma. Also, part (C) is a direct consequence of Corollary 3.4 and parts (A) and (B).

Finally, we consider parts (A) and (B) of the lemma. Fix $m \in \{1,\ldots,\mathbf{K}+\mathbf{J}-\mathbf{B}\}$. From (3.35), we have, on recalling that $U_m^r(\cdot)$ is nondecreasing and nonnegative, that for $0 \leq s \leq t$ and all $r$,

$$\tau^r(t) - \tau^r(s) \geq (t-s) \geq 0, \qquad \tau^r(t) - \tau^r(s) \geq \hat{U}_m^r(t) - \hat{U}_m^r(s) \geq 0,$$

and so

$$0 \leq \check{\tau}^r(t) - \check{\tau}^r(s) \leq (t-s), \qquad 0 \leq \check{U}_m^r(t) - \check{U}_m^r(s) \leq (t-s). \tag{3.39}$$

This proves that $\check{\tau}^r(\cdot)$ and $\check{U}^r(\cdot)$ have continuous sample paths and that $\{\check{\tau}^r(\cdot)\}_r$ and $\{\check{U}^r(\cdot)\}_r$ are $\mathcal{C}$-tight. The $\mathcal{C}$-tightness of $\{\check{X}^r(\cdot)\}_r$ and $\{\check{W}^r(\cdot)\}_r$ is now immediate from Lemma 3.3(C), on noting that $\check{X}^r(\cdot) = \hat{X}^r(\check{\tau}^r(\cdot))$ and $\check{W}^r(\cdot) = \Lambda\hat{q}^r + \Lambda\theta\check{\tau}^r(\cdot) + \Lambda\check{X}^r(\cdot) + G\check{U}^r(\cdot)$. Also, the $\mathcal{C}$-tightness of $\{\hat{M}^r, \check{M}^r\}$ follows from Corollary 3.4. This proves part (A) of the lemma.

For the proof of part (B) of the lemma, note that the continuity and nondecreasing property of $\check{\tau}$ are immediate consequences of the same properties of $\check{\tau}^r$, for each fixed $r$. It now suffices to show (3.37), since the remaining properties are then immediate. Observing that $\check{\tau}^r$ is nondecreasing and $\check{\tau}^r \Rightarrow$

CONTROLLED STOCHASTIC NETWORKS 25

$\check{\tau}$ as $r \to \infty$, it follows that in order to prove (3.37), it suffices to show for all $M > 0$,

$$\lim_{t \to \infty} \limsup_{r \to \infty} \mathbb{P}[\check{\tau}^r(t) < M] = 0. \tag{3.40}$$

Using Markov's inequality, we see that

$$\mathbb{P}[\check{\tau}^r(t) < M] = \mathbb{P}[\tau^r(M) > t] = \mathbb{P}\left[M + \sum_{m=1}^{(\mathbf{K}+\mathbf{J}-\mathbf{B})} \hat{U}_m^r(M) > t\right]$$

$$\leq \sum_{m=1}^{(\mathbf{K}+\mathbf{J}-\mathbf{B})} \mathbb{P}\left[\hat{U}_m^r(M) > \frac{(t-M)}{(\mathbf{K}+\mathbf{J}-\mathbf{B})}\right] \tag{3.41}$$

$$\leq \sum_{m=1}^{(\mathbf{K}+\mathbf{J}-\mathbf{B})} (\mathbf{K}+\mathbf{J}-\mathbf{B})\frac{\mathbb{E}[\hat{U}_m^r(M)]}{(t-M)}.$$

The statement in (3.40) now follows, using (3.21) and passing to the limit in (3.41), first as $r \to \infty$ and then as $t \to \infty$. $\square$

The following theorem is at the heart of our analysis. Since the proof is rather long, we postpone it to Section 4.

THEOREM 3.7. *Let* $(\check{W}, \check{X}, \check{U}, \check{\tau}, \tau, \Xi, M', \check{M})$, *given on the probability space* $(\hat{\Omega}, \hat{\mathcal{F}}, \hat{\mathbb{P}})$, *be as in Lemma* 3.6. *Define*

$$W(\cdot) \doteq \check{W}(\tau(\cdot)), \qquad X(\cdot) \doteq \check{X}(\tau(\cdot)),$$
$$U(\cdot) \doteq \check{U}(\tau(\cdot)), \qquad M(\cdot) \doteq \check{M}(\tau(\cdot)). \tag{3.42}$$

*Then, almost surely,* $M = M'$, *that is,*

$$M(\cdot) \equiv (M^E(\cdot), M^S(\cdot), M^\Phi(\cdot))$$
$$= \left(\tilde{E}(\cdot), i = 1, \ldots, \mathbf{I}, \tilde{S}_j(x_j^* \cdot), j = 1, \ldots, \mathbf{J}, \right. \tag{3.43}$$
$$\left. \sum_{j=1}^{\mathbf{J}} \tilde{\Phi}_i^j(\beta_j x_j^* \cdot), i = 1, \ldots, \mathbf{I}\right).$$

*Furthermore, there is a filtration* $\{\hat{\mathcal{F}}(t)\}_{t \geq 0}$ *on the probability space* $(\hat{\Omega}, \hat{\mathcal{F}}, \hat{\mathbb{P}})$ *such that the processes* $(W, X, U)$ *are adapted to* $\{\hat{\mathcal{F}}(t)\}$ *and* $X$ *is a* $\{\hat{\mathcal{F}}(t)\}$-*Brownian motion with zero drift and covariance matrix* $\Sigma$ *given by* (2.38). *That is,* $\hat{\Gamma} = (\hat{\Omega}, \hat{\mathcal{F}}, \hat{\mathbb{P}}, \{\hat{\mathcal{F}}(t)\}, X)$ *is a system in the sense of Definition* 2.9.

REMARK 3.8. From the relation $\check{X}^r(\cdot) = \hat{X}^r(\check{\tau}^r(\cdot))$, Lemma 3.3(C) and Corollary 3.4, it follows that $\check{X}(\cdot) = \tilde{X}(\check{\tau}(\cdot))$, a.s., where $\tilde{X}$ is defined via (3.7). Combining this with the fact that $\check{\tau}(\tau(t)) = t$, $t \geq 0$, we obtain that $X(\cdot) \doteq \check{X}(\tau(\cdot)) = \tilde{X}(\cdot)$. This immediately shows that $X$ is a $(0, \Sigma)$-Brownian motion and hence a martingale with respect to its own filtration. However, it is not clear that the limiting control process $U(\cdot)$ will be adapted to this filtration (i.e., the filtration generated by $\{X(t)\}$). Theorem 3.7 shows that there is a filtration to which all the limit processes are adapted, and with respect to which the limit $X$ is still a Brownian motion.

COROLLARY 3.9. *Let $\hat{\Gamma}$ be the system obtained in Theorem 3.7. Define $w = \Lambda q$. Let $(W, U)$ be as in (3.42). Then $U \in \tilde{\mathcal{A}}_0(w)$ and (2.46) in Definition 2.10 holds with $(\tilde{W}, \tilde{X}, \tilde{U}_0)$ replaced by $(W, X, U)$. In particular, we have*

$$(3.44) \qquad \hat{\mathbb{E}} \int_0^\infty e^{-\gamma t} \hat{h}(W(t)) \, dt + \hat{\mathbb{E}} \int_{[0,\infty)} e^{-\gamma t} p \cdot dU(t) \geq \tilde{J}_0^*(w),$$

*where $\hat{\mathbb{E}}$ denotes the expectation operator corresponding to $\hat{\mathbb{P}}$.*

PROOF. Since $\check{U}$ and $\tau$ are both nondecreasing, it follows that $U(\cdot) \doteq \check{U}(\tau(\cdot))$ is nondecreasing. Also, since $\hat{U}^r(t) \in \mathcal{K}$, for all $r$ and $t \geq 0$, we obtain that $U(t) \in \mathcal{K}$, $t \geq 0$. Using (2.24) and (3.36), we have that for all $t \geq 0$,

$$(3.45) \quad \begin{aligned} \check{W}^r(t) &= \Lambda \hat{q}^r + \Lambda[\theta_1^r \check{\tau}^r(t) - (C - P') \operatorname{diag}(\theta_2^r) \bar{T}^r(\check{\tau}^r(t))] \\ &\quad + \Lambda \check{X}^r(t) + G \check{U}^r(t) \in \mathcal{W}. \end{aligned}$$

Taking the limit as $r \to \infty$, one obtains that for all $t \geq 0$,

$$(3.46) \qquad \check{W}(t) = w + \Lambda \theta \check{\tau}(t) + \Lambda \check{X}(t) + G\check{U}(t) \in \mathcal{W}.$$

Thus, from (3.42), and recalling $\check{\tau}(\tau(t)) = t$, $t \geq 0$, we have for all $t \geq 0$,

$$(3.47) \qquad W(t) = w + \Lambda \theta t + \Lambda X(t) + GU(t) \in \mathcal{W}.$$

From Theorem 3.7, we have that $X(\cdot)$ is a $\{\hat{\mathcal{F}}(t)\}$-Brownian motion and $U(\cdot), W(\cdot)$ are $\{\hat{\mathcal{F}}(t)\}$-adapted. Hence, $U \in \tilde{\mathcal{A}}_0(w)$, with associated system $\hat{\Gamma}$. The inequality in (3.44) is now an immediate consequence of the definition of the value function of the EWF, $\tilde{J}_0^*(w)$, in (2.48). □

For the proof of the following lemma, we refer the reader to Theorem IV.4.5 of [21].

LEMMA 3.10. *Let $a$ be a $\mathbb{R}_+$-valued, right-continuous function on $[0, \infty)$ such that $a(0) = 0$. Let $c$ be its right inverse, that is, $c(t) \doteq \inf\{s \geq 0 : a(s) > t\}$, $t \geq 0$. Assume that $c(t) < \infty$ for all $t \geq 0$. Let $f$ be a nonnegative Borel* *measurable function on* $[0,\infty)$. *Then if* $G$ *is any* $\mathbb{R}_+$-*valued, right-continuous function on* $[0,\infty)$,

$$(3.48) \quad \int_{[0,\infty)} f(s)\, dG(a(s)) = \int_{[0,\infty)} f(c(s-))\, dG(s),$$

*with the convention that the contribution to the integrals above at* $0$ *is* $f(0)G(0)$. *This, in particular, implies*

$$(3.49) \quad \int_{[0,\infty)} f(s)\, da(s) = \int_{[0,\infty)} f(c(s))\, ds.$$

The following lemma will be needed in our proofs. We refer the reader to [8] (Lemma 2.4) for a proof. Let $\mathcal{C}_b^d$ denote the space of all bounded functions in $\mathcal{C}^d$.

LEMMA 3.11. *Suppose* $\xi^r \to \xi$ *in* $\mathcal{D}^d$ *and* $\lambda^r \to \lambda$ *in* $\mathcal{C}$ *as* $r \to \infty$. *Further, suppose that* $\lambda^r$ *is nonnegative and nondecreasing for each* $r$. *Then for any* $f \in \mathcal{C}_b^d =$ *space of bounded functions in* $\mathcal{C}^d$,

$$\int_{[0,u)} f(\xi^r(t))\, d\lambda^r(t) \to \int_{[0,u)} f(\xi(t))\, d\lambda(t) \qquad as\ r \to \infty,$$

*uniformly for all* $u$ *in any compact subset of* $[0,\infty)$.

PROOF OF THEOREM 3.1. On noting that $\check{\tau}^r(\cdot)$ is continuous and $\check{\tau}^r(t) \uparrow \infty$ as $t \to \infty$, we have that

$$(3.50) \quad \begin{aligned} J^r(q^r, T^r) &= \mathbb{E} \int_0^\infty e^{-\gamma t} h \cdot \hat{Q}^r(t)\, dt + \mathbb{E}\left(\int_0^\infty e^{-\gamma t} p \cdot d\hat{U}^r(t)\right) \\ &= \mathbb{E} \int_0^\infty e^{-\gamma \check{\tau}^r(t)} h \cdot \check{Q}^r(t)\, d\check{\tau}^r(t) + \mathbb{E}\left(\int_0^\infty e^{-\gamma \check{\tau}^r(t)} p \cdot d\check{U}^r(t)\right), \end{aligned}$$

where the last equality follows from Lemma 3.10, with $a = \check{\tau}^r$ and $c = \tau^r$ [(3.49) is used for the first term and (3.48) with $G = \hat{U}^r$ is used for the second term]. Recall the definitions of $\hat{h}$ and $w(\cdot)$ from Section 2. Define

$$(3.51) \quad f(x,y) \doteq e^{-\gamma y} \hat{h}(x), \qquad x \in \mathbb{R}_+^{\mathbf{L}}, y \in \mathbb{R}_+.$$

Since $h \cdot q \geq \hat{h}(w(q))$ for all $q \in \mathbb{R}_+^{\mathbf{I}}$, we have that the first term right-hand side of (3.50) is bounded from below by

$$(3.52) \quad \mathbb{E} \int_0^\infty e^{-\gamma \check{\tau}^r(t)} \hat{h}(\check{W}^r(t))\, d\check{\tau}^r(t) = \mathbb{E} \int_0^\infty f(\check{W}^r(t), \check{\tau}^r(t))\, d\check{\tau}^r(t).$$

Let $(\check{W}(\cdot), \check{X}(\cdot), \check{U}(\cdot), \check{\tau}(\cdot), \Xi(\cdot), M'(\cdot), \check{M}(\cdot))$ be as in Lemma 3.6 and let $\{r'\}$ be the sequence along which convergence to this limit point occurs. Using the Skorohod representation theorem (and relabeling the subsequence $\{r'\}$ by $\{r\}$), we can assume without loss of generality that

$$
\begin{aligned}
(\check{W}^r(\cdot), \check{U}^r(\cdot)) &\to (\check{W}(\cdot), \check{U}(\cdot)) \quad \text{in } \mathcal{D} \quad \text{and} \\
\check{\tau}^r(\cdot) &\to \check{\tau}(\cdot) \quad \text{in } \mathcal{C} \text{ a.s. as } r \to \infty.
\end{aligned}
\tag{3.53}
$$

Define, for $N \geq 1$, $f_N(\cdot,\cdot) \doteq f(\cdot,\cdot) \wedge N$. Using Lemma 3.11, we have for each fixed $N \geq 1$ and $u \geq 0$,

$$
\int_0^u f_N(\check{W}^r(t), \check{\tau}^r(t))\, d\check{\tau}^r(t) \to \int_0^u f_N(\check{W}(t), \check{\tau}(t))\, d\check{\tau}(t) \quad \text{a.s.,}
\tag{3.54}
$$

$$
\int_0^u e^{-\gamma \check{\tau}^r(t)} p \cdot d\check{U}^r(t) \to \int_0^u e^{-\gamma \check{\tau}(t)} p \cdot d\check{U}(t) \quad \text{a.s.}
\tag{3.55}
$$

as $r \to \infty$. Thus, a.s.,

$$
\begin{aligned}
\liminf_{r \to \infty} \int_0^\infty f(\check{W}^r(t), \check{\tau}^r(t))\, d\check{\tau}^r(t) \\
\geq \liminf_{r \to \infty} \int_0^u f_N(\check{W}^r(t), \check{\tau}^r(t))\, d\check{\tau}^r(t) \\
= \int_0^u f_N(\check{W}(t), \check{\tau}(t))\, d\check{\tau}(t).
\end{aligned}
\tag{3.56}
$$

Taking the limit as $N \to \infty$ and $u \to \infty$ in (3.56), we have that, a.s.,

$$
\liminf_{r \to \infty} \int_0^\infty f(\check{W}^r(t), \check{\tau}^r(t))\, d\check{\tau}^r(t) \geq \int_0^\infty f(\check{W}(t), \check{\tau}(t))\, d\check{\tau}(t).
\tag{3.57}
$$

Similarly, using (3.55), we obtain

$$
\liminf_{r \to \infty} \int_0^\infty e^{-\gamma \check{\tau}^r(t)} p \cdot d\check{U}^r(t) \geq \int_0^\infty e^{-\gamma \check{\tau}(t)} p \cdot d\check{U}(t).
\tag{3.58}
$$

Next, using (3.4), (3.50) and (3.52), we obtain

$$
\begin{aligned}
J(q, \{T^r\}) &= \lim_{r \to \infty} J^r(q^r, T^r) \\
&= \lim_{r \to \infty} \mathbb{E} \int_0^\infty e^{-\gamma \check{\tau}^r(t)} h(\check{Q}^r(t))\, d\check{\tau}^r(t) \\
&\quad + \lim_{r \to \infty} \mathbb{E} \int_0^\infty e^{-\gamma \check{\tau}^r(t)} p \cdot d\check{U}^r(t) \\
&\geq \liminf_{r \to \infty} \mathbb{E} \int_0^\infty f(\check{W}^r(t), \check{\tau}^r(t))\, d\check{\tau}^r(t) \\
&\quad + \liminf_{r \to \infty} \mathbb{E} \int_0^\infty e^{-\gamma \check{\tau}^r(t)} p \cdot d\check{U}^r(t).
\end{aligned}
\tag{3.59}
$$

Finally, using Fatou's lemma, (3.57), (3.58), (3.59) and (3.51), we have

$$
\begin{aligned}
J(q, \{T^r\}) &\geq \mathbb{E}\bigg[\liminf_{r\to\infty}\int_0^\infty f(\check{W}^r(t), \check{\tau}^r(t))\,d\check{\tau}^r(t) \\
&\qquad + \liminf_{r\to\infty}\int_0^\infty e^{-\gamma\check{\tau}^r(t)}p\cdot d\check{U}^r(t)\bigg] \\
&\geq \hat{\mathbb{E}}\bigg[\int_0^\infty f(\check{W}(t), \check{\tau}(t))\,d\check{\tau}(t) + \int_0^\infty e^{-\gamma\check{\tau}(t)}p\cdot d\check{U}(t)\bigg] \\
&= \hat{\mathbb{E}}\bigg[\int_0^\infty f(W(t), t)\,dt + \int_{[0,\infty)} e^{-\gamma t}p\cdot dU(t)\bigg] \\
&= \hat{\mathbb{E}}\bigg[\int_0^\infty e^{-\gamma t}\hat{h}(W(t))\,dt + \int_{[0,\infty)} e^{-\gamma t}p\cdot dU(t)\bigg] \\
&\geq \tilde{J}_0^*(w),
\end{aligned}
$$
(3.60)

where the equalities in the third line follow from Lemma 3.10 [with $a = \check{\tau}$, $c = \tau$ in (3.49) for the first term and $a = \tau, c = \check{\tau}$, $G = \check{U}$ in (3.48) for the second term]. The inequality in the last line of (3.60) is a consequence of Corollary 3.9. This completes the proof. □

**4. Proof of Theorem 3.7.** This section is devoted to the proof of Theorem 3.7. Since $\hat{U}^r(\cdot)$ is adapted to $\{\mathcal{F}_1^r(t) : t \geq 0\}$, it follows that for all $s, a \geq 0$,

$$\{\check{\tau}^r(s) \leq a\} = \{\tau^r(a) \geq s\} = \bigg\{a + \sum_{m=1}^{\mathbf{K}+\mathbf{J}-\mathbf{B}} \hat{U}_m^r(a) \geq s\bigg\} \in \mathcal{F}_1^r(a).$$

This shows that for each $s \geq 0$, $\check{\tau}^r(s)$ is a $\{\mathcal{F}_1^r(t) : t \geq 0\}$ stopping time. For each $r, t \geq 0$, define the stopped sigma field

(4.1) $$\mathcal{G}^r(t) \doteq \mathcal{F}_1^r(\check{\tau}^r(t)).$$

Since, for each $r$, $\{\check{\tau}^r(t) : t \geq 0\}$ is an increasing sequence of stopping times, it follows that $\{\mathcal{G}^r(t) : t \geq 0\}$ is a filtration. The following lemma gives an alternative representation of the filtration $\{\mathcal{G}^r(t)\}$:

LEMMA 4.1. *Define for fixed $r$ and $t \geq 0$, $\rho^r(t) \doteq \sigma_0^r(\check{\tau}^r(t))$. Then $\rho^r(t)$ is a $\{\mathcal{F}^r((m,n))\}$ stopping time,*

(4.2) $$\mathcal{G}^r(t) = \mathcal{F}^r(\rho^r(t))$$

*and $\check{U}^r(\cdot) \doteq \hat{U}^r(\check{\tau}^r(\cdot))$ is $\{\mathcal{G}^r(t)\}$-adapted.*

PROOF. Fix $r, t \geq 0$. Let $\{\check{\tau}_\alpha^r(t)\}_\alpha$ be a decreasing sequence of $\{\mathcal{F}_1^r(s)\}$ stopping times, such that for each $\alpha \in \mathbb{N}$, $\check{\tau}_\alpha^r(t)$ takes values in the countable set $\Gamma_\alpha \doteq \{\frac{j}{2^\alpha} : j \in \mathbb{N}_0\}$ and $\check{\tau}^r_\alpha(t) \downarrow \check{\tau}^r(t)$ as $\alpha \to \infty$. By right continuity of $\sigma^r_0(\cdot)$, it follows that $\rho^r_\alpha(t) \doteq \sigma^r_0(\check{\tau}^r_\alpha(t)) \downarrow \rho^r(t)$, as $\alpha \to \infty$.

Next, observe that for all $(m, n) \in \mathbb{N}^{\mathbf{I+J}}$, we can write

$$\{\rho^r_\alpha(t) \leq (m, n)\} = \bigcup_{s \in \Gamma_\alpha} \{\sigma^r_0(s) \leq (m, n)\} \cap \{\check{\tau}^r_\alpha(t) = s\}. \tag{4.3}$$

For each $s \in \Gamma_\alpha$, $\{\check{\tau}^r_\alpha(t) = s\} \in \mathcal{F}^r_1(s) = \mathcal{F}^r(\sigma^r_0(s))$. Thus, by definition of the stopped sigma field $\mathcal{F}^r(\sigma^r_0(s))$, the set in (4.3) is in $\mathcal{F}^r((m,n))$ for all $m, n$, and so for each $\alpha$, $\rho^r_\alpha(t)$ is a $\{\mathcal{F}^r((m,n))\}$ stopping time. Thus, $\rho^r(t) = \inf_\alpha \rho^r_\alpha(t)$ is a $\{\mathcal{F}^r((m,n))\}$ stopping time.

To prove (4.2), define $\tilde{\mathcal{G}}^r(t) \doteq \mathcal{F}^r(\rho^r(t))$ and for each $\alpha$, $\mathcal{G}^r_\alpha(t) \doteq \mathcal{F}^r_1(\check{\tau}^r_\alpha(t))$ and $\tilde{\mathcal{G}}^r_\alpha(t) \doteq \mathcal{F}^r(\rho^r_\alpha(t))$ for $t \geq 0$. Next, we show that for $r, \alpha, t \geq 0$,

$$\tilde{\mathcal{G}}^r_\alpha(t) = \mathcal{G}^r_\alpha(t). \tag{4.4}$$

Consider $A \in \tilde{\mathcal{G}}^r_\alpha(t)$. Then

$$A \cap \{\rho^r_\alpha(t) = (m, n)\} \in \mathcal{F}^r((m, n)) \qquad \text{for all } (m, n) \in \mathbb{N}^{\mathbf{I+J}}. \tag{4.5}$$

Since $\check{\tau}^r_\alpha(t)$ is a $\{\mathcal{F}^r_1(s) : s \geq 0\}$ stopping time taking values in $\Gamma_\alpha$, $\{\check{\tau}^r_\alpha(t) = s\} \in \mathcal{F}^r_1(s) \doteq \mathcal{F}^r(\sigma^r_0(s))$ for all $s \in \Gamma_\alpha$. Hence, for all $s \in \Gamma_\alpha, (m, n) \in \mathbb{N}^{\mathbf{I+J}}$,

$$\{\check{\tau}^r_\alpha(t) = s\} \cap \{\sigma^r_0(s) = (m, n)\} \in \mathcal{F}^r((m, n)). \tag{4.6}$$

Taking the intersection of the two sets in (4.5) and (4.6), and recalling the definition of $\rho^r_\alpha(t)$, we obtain that for all $s \in \Gamma_\alpha, (m, n) \in \mathbb{N}^{\mathbf{I+J}}$,

$$A \cap \{\check{\tau}^r_\alpha(t) = s\} \cap \{\sigma^r_0(s) = (m, n)\} \in \mathcal{F}^r((m, n)). \tag{4.7}$$

Thus, $A \cap \{\check{\tau}^r_\alpha(t) = s\} \in \mathcal{F}^r_1(s)$ for all $s \in \Gamma_\alpha$, proving $A \in \mathcal{G}^r_\alpha(t)$. This proves $\tilde{\mathcal{G}}^r_\alpha(t) \subseteq \mathcal{G}^r_\alpha(t)$. Next, we consider the reverse inclusion. If $A \in \mathcal{G}^r_\alpha(t)$, then (4.7) holds, and taking the union over all $s \in \Gamma_\alpha$, one gets (4.5), which implies that $A \in \tilde{\mathcal{G}}^r_\alpha(t)$. This completes the proof of (4.4).

Note that if $\gamma_k$ is a $\{\mathcal{F}^r((m,n))\}$ stopping time with $\gamma_k \downarrow \gamma$, then $\gamma$ is a $\{\mathcal{F}^r((m,n))\}$ stopping time and $\mathcal{F}^r(\gamma_k) \downarrow \mathcal{F}^r(\gamma)$. Using this property and noting that $\sigma^r_0(t)$ is a right-continuous function of $t$, it follows that $\{\mathcal{F}^r_1(t) = \mathcal{F}^r(\sigma^r_0(t)) : t \geq 0\}$ is a right-continuous filtration. This, together with (4.4) and the facts that $\check{\tau}^r_\alpha(t) \downarrow \check{\tau}^r(t)$ and $\rho^r_\alpha(t) \downarrow \rho^r(t)$, yields (4.2). The adaptedness of $\check{U}^r(\cdot)$ to $\{\tilde{\mathcal{G}}^r(t)\}$ is an immediate consequence of the adaptedness of $\hat{U}^r$ to $\{\mathcal{F}^r_1(t)\}$ and (4.1) (cf. Proposition 1.2.18 of [14]). □

The following lemma is needed for some estimates in the proof of Theorem 3.7:

LEMMA 4.2. *For each $t \in [0, \infty)$ and $p \in \mathbb{N}$,*

$$\sup_r \max_{i=1,2,\ldots,\mathbf{I}} \mathbb{E}((E^r_i(t))^p) < \infty, \qquad \sup_r \max_{j=1,2,\ldots,\mathbf{J}} \mathbb{E}((S^r_j(t))^p) < \infty.$$

PROOF. We will only prove the first inequality. Fix $i \in \{1, 2, \ldots, \mathbf{I}\}$. We claim that there exist $\varepsilon, \delta \in (0, 1)$ such that

(4.8) $$\inf_r \mathbb{P}(u_i^r(1) > \delta) > \varepsilon.$$

To prove the claim, we argue via contradiction. In view of Assumption 2.1, there exists an $a \in (0, \infty)$ such that

(4.9) $$\inf_r \mathbb{E}(u_i^r(1)) \geq a.$$

Let $\delta \doteq (a/4) \wedge (1/2)$ and suppose that (4.8) fails to hold for this choice of $\delta$, for any $\varepsilon \in (0, 1)$. Then there exists a sequence $\{r_k\}_{k \geq 1}$ such that

(4.10) $$\mathbb{P}(u_i^{r_k}(1) > \delta) \leq \frac{1}{k}.$$

This, in particular, implies that $i \leq \mathbf{I}'$ and that $\{u_i^{r_k}(1)\}_{k \geq 1}$ is a tight sequence. In fact, since $\sup_r \mathbb{E}(u_i^r(1))^2 < \infty$ from Assumption 2.1, the sequence is uniformly integrable. Let $u^*$ be a limit point of this sequence along a weakly convergent subsequence. By uniform integrability and (4.9), we have $\mathbb{E}(u^*) \geq a$. Thus, $\mathbb{P}(u^* > 2\delta) > \varepsilon$ for some $\varepsilon \in (0, \infty)$. However, in view of (4.10), and recalling that $u^*$ is a (weak) limit point of $\{u_i^{r_k}(1)\}_{k \geq 1}$, we see that $\mathbb{P}(u^* > \delta) = 0$. Thus, we arrive at a contradiction. This proves (4.8). Now choose $\varpi \in \mathbb{N}$ such that $\varpi \delta \geq t$. Following the proof of Theorem 3.4.2 of [9], we see that for all $m \in \mathbb{N}_0$ and all $r$, $\mathbb{P}(E_i^r(t) > m\varpi) \leq (1 - \varepsilon^\varpi)^m$. The moment estimate in the statement of the lemma now follows. $\square$

*Martingales associated with $\check{M}^r(\cdot)$.* Define for $m \in \mathbb{N}^{\mathbf{I}}$, $n \in \mathbb{N}^{\mathbf{J}}$, $i' = 1, \ldots, \mathbf{I}'$, $j = 1, \ldots, \mathbf{J}$, $i = 1, \ldots, \mathbf{I}$,

(4.11)
$$\tilde{\xi}_{i'}^r(m, n) \doteq \tilde{\xi}_{i'}^r(m_{i'}) = \frac{1}{r} \sum_{n=1}^{m_{i'}} (1 - \alpha_{i'}^r u_i^r(n)), \qquad m_{i'} \geq 1,$$

$$\tilde{\eta}_j^r((m, n)) \doteq \tilde{\eta}_j^r(n_j) = \frac{1}{r} \sum_{n=1}^{n_j} (1 - \beta_j^r v_j^r(n)), \qquad n_j \geq 1,$$

$$\tilde{\zeta}_i^r((m, n)) \doteq \sum_{j=1}^{\mathbf{J}} \tilde{\zeta}_i^{j,r}(n_j),$$

$$\tilde{\zeta}_i^{j,r}(n_j) \doteq \frac{1}{r} \sum_{k=1}^{n_j} (\phi_i^{j,r}(k) - p_i^j), \qquad n_j \geq 1.$$

As a convention, we take $\tilde{\xi}_i^r((m, n)) = 0$ for $i = \mathbf{I}' + 1, \ldots, \mathbf{I}$. It is easy to check that for all $i = 1, \ldots, \mathbf{I}$, $j = 1, \ldots, \mathbf{J}$, $\tilde{\xi}_i^r(\cdot)$, $\tilde{\eta}_j^r(\cdot)$ and $\tilde{\zeta}_i^r(\cdot)$ are $\{\mathcal{F}^r((m, n)) : m \in \mathbb{N}^{\mathbf{I}}, n \in \mathbb{N}^{\mathbf{J}}\}$ martingales. Straightforward calculations give the following quadratic variations for these martingales. For $i, i_1 \neq i_2 \in \{1, \ldots, \mathbf{I}\}, j, j_1 \neq j_2 \in \{1, \ldots, \mathbf{J}\}$,

$$\langle \tilde{\xi}_i^r \rangle((m,n)) = \frac{m_i(\alpha_i^r \sigma_i^{u,r})^2}{r^2},$$

$$\langle \tilde{\eta}_j^r \rangle((m,n)) = \frac{n_j(\beta_j^r \sigma_j^{v,r})^2}{r^2},$$

$$\langle \tilde{\zeta}_i^r \rangle((m,n)) = \frac{\sum_{j=1}^{\mathbf{J}} n_j(\sigma_{ii}^{\phi^j})}{r^2},$$

(4.12)

$$\langle \tilde{\xi}_{i_1}^r, \tilde{\xi}_{i_2}^r \rangle((m,n)) = \langle \tilde{\xi}_i^r, \tilde{\eta}_j^r \rangle((m,n)) = 0,$$

$$\langle \tilde{\eta}_{j_1}^r, \tilde{\eta}_{j_2}^r \rangle((m,n)) = \langle \tilde{\eta}_j^r, \tilde{\zeta}_{i_1}^r \rangle((m,n)) = \langle \tilde{\xi}_i^r, \tilde{\zeta}_{i_1}^r \rangle((m,n)) = 0,$$

$$\langle \tilde{\zeta}_{i_1}^r, \tilde{\zeta}_{i_2}^r \rangle((m,n)) = \frac{\sum_{j=1}^{\mathbf{J}} n_j(\sigma_{i_1 i_2}^{\phi^j})}{r^2}.$$

Let $\rho^r(t)$ be as in Lemma 4.1. Using the fact that $\rho^r(t)$ is a $\{\mathcal{F}^r((m,n))\}$ (multiparameter) stopping time, we have that

(4.13)
$$\check{N}^r(t) = (\check{N}_i^{E,r}(t), i = 1, \ldots, \mathbf{I};$$
$$\check{N}_j^{S,r}(t), j = 1, \ldots, \mathbf{J}; \check{N}_i^{\Phi,r}(t), i = 1, \ldots, \mathbf{I})$$
$$\doteq (\tilde{\xi}_i^r(\rho^r(t)), i = 1, \ldots, \mathbf{I};$$
$$\tilde{\eta}_j^r(\rho^r(t)), j = 1, \ldots, \mathbf{J}; \tilde{\zeta}_i^r(\rho^r(t)), i = 1, \ldots, \mathbf{I})$$

is a $\{\mathcal{G}^r(t) : t \geq 0\}$ martingale. This multiparameter version of optional sampling theorem can be proven in a manner similar to the corresponding single parameter case (cf. Theorems 4.7.4 and 4.7.5 of [9]). The key conditions needed to invoke the theorem are the following. For all $i = 1, \ldots, \mathbf{I}$, $j = 1, \ldots, \mathbf{J}$, $m_i, n_j = 1, 2, \ldots$,

(4.14)
$$\mathbb{E}[|\rho^r(t)|] < \infty,$$
$$\mathbb{E}[|\tilde{\xi}_i^r(m_i+1) - \tilde{\xi}_i^r(m_i)||\mathcal{F}^r((m,n))] \leq c^r,$$
$$\mathbb{E}[|\tilde{\eta}_j^r(n_j+1) - \tilde{\eta}_j^r(n_j)||\mathcal{F}^r((m,n))] \leq c^r,$$
$$\mathbb{E}[|\tilde{\zeta}_i^{j,r}(n_j+1) - \tilde{\zeta}_i^{j,r}(n_j)||\mathcal{F}^r((m,n))] \leq c^r,$$

where $c^r \in (0, \infty)$ does not depend on $(m,n)$. The last three bounds in (4.14) follow from the definition of the martingales and the fact that for fixed $r$, $i$ and $j$, $\{u_i^r(n), v_j^r(n), \phi_i^{j,r}(n)\}$ is a sequence of i.i.d. random variables with finite variance. For the first condition in (4.14), note that using $\check{\tau}^r(t) \leq t$,

$\bar{T}_j^r(t) \le t$ for all $t \ge 0$, we have, for some $\kappa \in (0, \infty)$,

$$\mathbb{E}[|\rho^r(t)|] \le \kappa \left[ \sum_{i=1}^{\mathbf{I}} \mathbb{E}(E_i^r(r^2 \tilde{\tau}^r(t) + 1)) + \sum_{j=1}^{\mathbf{J}} \mathbb{E}(S_j^r(T_j^r(r^2 \tilde{\tau}^r(t)) + 1)) \right]$$

$$(4.15) \qquad \le \kappa \left[ \sum_{i=1}^{\mathbf{I}} \mathbb{E}(E_i^r((r^2 t) + 1)) + \sum_{j=1}^{\mathbf{J}} \mathbb{E}(S_j^r((r^2 t) + 1)) \right]$$

$$< \infty.$$

This proves (4.14) and hence the martingale property claimed below (4.13).

Using (4.12) and the above martingale property, it is easy to verify that the following hold for any $t_2^r \ge t_1^r$ and all $i, i_1, i_2, \in \{1, \ldots, \mathbf{I}\}$, $j \in \{1, \ldots, \mathbf{J}\}$:

$$(4.16) \quad \begin{aligned} & \mathbb{E}[(\check{N}_i^{E,r}(t_2^r) - \check{N}_i^{E,r}(t_1^r))^2 | \mathcal{G}^r(t_1^r)] \\ &= (\alpha_i^r \sigma_i^{u,r})^2 \mathbb{E}[(\bar{E}_i^r(\tilde{\tau}^r(t_2^r)) - \bar{E}_i^r(\tilde{\tau}^r(t_1^r))) | \mathcal{G}^r(t_1^r)], \end{aligned}$$

$$(4.17) \quad \begin{aligned} & \mathbb{E}[(\check{N}_j^{S,r}(t_2^r) - \check{N}_j^{S,r}(t_1^r))^2 | \mathcal{G}^r(t_1^r)] \\ &= (\beta_j^r \sigma_i^{v,r})^2 \mathbb{E}[(\bar{S}_j^r(\bar{T}_j^r(\tilde{\tau}^r(t_2^r))) - \bar{S}_j^r(\bar{T}_j^r(\tilde{\tau}^r(t_1^r)))) | \mathcal{G}^r(t_1^r)], \end{aligned}$$

$$(4.18) \quad \begin{aligned} & \mathbb{E}[(\check{N}_{i_1}^{\Phi,r}(t_2^r) - \check{N}_{i_1}^{\Phi,r}(t_1^r))(\check{N}_{i_2}^{\Phi,r}(t_2^r) - \check{N}_{i_2}^{\Phi,r}(t_1^r)) | \mathcal{G}^r(t_1^r)] \\ &= \sum_{j=1}^{\mathbf{J}} (\sigma_{i_1 i_2}^{\phi^j}) \mathbb{E}[(\bar{S}_j^r(\bar{T}_j^r(\tilde{\tau}^r(t_2^r))) - \bar{S}_j^r(\bar{T}_j^r(\tilde{\tau}^r(t_1^r)))) | \mathcal{G}^r(t_1^r)]. \end{aligned}$$

Furthermore, for $i, i_1 \ne i_2 \in \{1, \ldots, \mathbf{I}\}$, $j_1 \ne j_2, \in \{1, \ldots, \mathbf{J}\}$,

$$(4.19) \quad \begin{aligned} & \mathbb{E}[(\check{N}_{i_1}^{E,r}(t_2^r) - \check{N}_{i_1}^{E,r}(t_1^r))(\check{N}_{i_2}^{E,r}(t_2^r) - \check{N}_{i_2}^{E,r}(t_1^r)) | \mathcal{G}^r(t_1^r)] \\ &= \mathbb{E}[(\check{N}_{j_1}^{S,r}(t_2^r) - \check{N}_{j_1}^{S,r}(t_1^r))(\check{N}_{j_2}^{S,r}(t_2^r) - \check{N}_{j_2}^{S,r}(t_1^r)) | \mathcal{G}^r(t_1^r)] \\ &= \mathbb{E}[(\check{N}_{i_1}^{E,r}(t_2^r) - \check{N}_{i_1}^{E,r}(t_1^r))(\check{N}_{j_1}^{S,r}(t_2^r) - \check{N}_{j_1}^{S,r}(t_1^r)) | \mathcal{G}^r(t_1^r)] \\ &= \mathbb{E}[(\check{N}_{i_1}^{\Phi,r}(t_2^r) - \check{N}_{i_1}^{\Phi,r}(t_1^r))(\check{N}_{j_1}^{S,r}(t_2^r) - \check{N}_{j_1}^{S,r}(t_1^r)) | \mathcal{G}^r(t_1^r)] \\ &= \mathbb{E}[(\check{N}_{i_1}^{\Phi,r}(t_2^r) - \check{N}_{i_1}^{\Phi,r}(t_1^r))(\check{N}_i^{E,r}(t_2^r) - \check{N}_i^{E,r}(t_1^r)) | \mathcal{G}^r(t_1^r)] \\ &= 0. \end{aligned}$$

*Relationship between $\check{M}^r(\cdot)$ and $\check{N}^r(\cdot)$.* Recall the definitions of $\check{M}^r$ and $\check{N}^r$ given in (3.36) [see also (3.17)] and (4.13), respectively. We now show that the difference between $\check{M}^r(\cdot)$ and $\check{N}^r(\cdot)$ approaches 0 as $r$ increases. Note that for $i = \mathbf{I}' + 1, \ldots, \mathbf{I}$, $\check{M}_i^{E,r}(\cdot) = \check{N}_i^{E,r}(\cdot) = 0$. Now consider $i \in \{1, \ldots, \mathbf{I}'\}$. It is easy to establish that for all $r$ and $t \ge 0$,

$$(4.20) \qquad \hat{E}_i^r(t) + \frac{1 - \alpha_i^r u_i^r(\sigma_{0,i}^{r,E}(t))}{r} \le \tilde{\xi}_i^r(\sigma_0^r(t)) \le \hat{E}_i^r(t) + \frac{1}{r}.$$

Recalling $\check{\tau}^r(t) \leq t$ for all $t \geq 0$, we see that

$$\sup_{0 \leq s \leq t} |\check{M}_i^{E,r}(s) - \check{N}_i^{E,r}(s)|$$
(4.21)
$$\leq \frac{1}{r} + \alpha_i^r \frac{[\sup_{0 \leq s \leq t}(u_i^r(E_i^r(r^2 s) + 1) - 1/\alpha_i^r)]}{r}.$$

Observing that $\{u_i^r(k) - 1/\alpha_i^r : k = 1, 2, \ldots\}$ is a sequence of i.i.d. zero mean random variables, we see that $\max_{k \leq (\alpha_i t+1)r^2} |u_i^r(k) - 1/\alpha_i^r|/r \to 0$ in probability, for all fixed $t \geq 0$. Using this observation, together with the fact that

(4.22) $$\frac{E_i^r(r^2 \cdot) + 1}{r^2} \Rightarrow \alpha_i e(\cdot) \qquad \text{as } r \to \infty,$$

it follows that

(4.23) $$\frac{[\sup_{0 \leq s \leq t}(u_i^r(E_i^r(r^2 s) + 1) - 1/\alpha_i^r)]}{r} \to 0 \qquad \text{as } r \to \infty,$$

in probability. Hence, using (4.21) and (4.23), we have for all $i = 1, \ldots, \mathbf{I}$,

(4.24) $$|\check{M}_i^{E,r}(\cdot) - \check{N}_i^{E,r}(\cdot)| \to 0 \qquad \text{as } r \to \infty,$$

in probability u.o.c.

Similar arguments, using the observations that $T_j^r(s) \leq s$ for all $s \geq 0$ and $\phi_i^{j,r}(n_{ij}), p_i^j \leq 1$ for all $n_{ij} \geq 1$, give the following two bounds:

$$\sup_{0 \leq s \leq t} |\check{M}_j^{S,r}(s) - \check{N}_j^{S,r}(s)|$$
(4.25)
$$\leq \frac{1}{r} + \beta_j^r \frac{[\sup_{0 \leq s \leq t}(v_j(S_j(r^2 s) + 1) - 1/\beta_j^r)]}{r},$$

$$\sup_{0 \leq s \leq t} |\check{M}_i^{\Phi,r}(s) - \check{N}_i^{\Phi,r}(s)| \leq \frac{2\mathbf{J}}{r}.$$

One can now argue, as before, that for $j = 1, \ldots, \mathbf{J}$, $i = 1, \ldots, \mathbf{I}$,

(4.26) $\quad |\check{M}_j^{S,r}(\cdot) - \check{N}_j^{S,r}(\cdot)| \to 0, \qquad |\check{M}_i^{\Phi,r}(\cdot) - \check{N}_i^{\Phi,r}(\cdot)| \to 0 \qquad \text{as } r \to \infty,$

in probability u.o.c. Thus, we have shown that

(4.27) $\quad |\check{M}^r(\cdot) - \check{N}^r(\cdot)| \to 0 \qquad \text{in probability u.o.c. as } r \to \infty.$

PROOF OF THEOREM 3.7. The first part of the theorem is immediate from Lemma 3.6(C) and the fact that $\check{\tau}(\tau(t)) = t$, $t \geq 0$. Using (2.21)–(2.24) and (3.18), we get, for all $t \geq 0$,

(4.28) $$\hat{X}^r(t) = \hat{M}^{E,r}(t) + (C - P')\hat{M}^{S,r}(t) + \hat{M}^{\Phi,r}(t),$$

(4.29) $$\hat{Y}^r(t) = r(x^* t - \bar{T}^r(t)), \qquad \hat{U}^r(t) = K\hat{Y}^r(t),$$

(4.30) $$\hat{W}^r(t) = \Lambda \hat{q}^r + \Lambda[\theta_1^r t - (C - P')\operatorname{diag}(\theta_2^r)\bar{T}^r(t)]$$
$$+ \Lambda \hat{X}^r(t) + G\hat{U}^r(t).$$

Replacing $t$ by $\check{\tau}^r(t)$ in the above, we get

(4.31) $$\check{X}^r(t) = \check{M}^{E,r}(t) + (C - P')\check{M}^{S,r}(t) + \check{M}^{\Phi,r}(t),$$

(4.32) $$\check{W}^r(t) = \Lambda \hat{q}^r + \Lambda[\theta_1^r \check{\tau}^r(t) - (C - P')\operatorname{diag}(\theta_2^r)\bar{T}^r(\check{\tau}^r(t))]$$
$$+ \Lambda \check{X}^r(t) + G\check{U}^r(t).$$

On passing to the limit in (4.31), (4.32) and using (3.42), we have

(4.33) $$X(t) = M^E(t) + (C - P')M^S(t) + M^\Phi(t),$$
$$W(t) = \Lambda q + \Lambda \theta t + \Lambda X(t) + GU(t),$$

for all $t \geq 0$, where $\theta$ is as in (2.39).

Define the sigma fields

(4.34) $$\check{\mathcal{H}}'(t) \doteq \sigma\{\check{M}(s), \check{U}(s), \check{\tau}(s) : s \leq t\},$$
$$\check{\mathcal{H}}(t) \doteq \bigcap_{n \geq 1} \check{\mathcal{H}}'\left(t + \frac{1}{n}\right), \qquad t \geq 0.$$

By construction, $\{\check{\mathcal{H}}(t)\}$ is a right-continuous filtration, and for all $u \geq 0$, $t \geq 0$, the following holds:

(4.35) $$\{\tau(t) < u\} = \{\check{\tau}(u) > t\} \in \check{\mathcal{H}}(u).$$

Using the right continuity of $\{\check{\mathcal{H}}(t)\}$, we now have that for all $t \geq 0$, $\tau(t)$ is a $\{\check{\mathcal{H}}(t)\}$ stopping time. Define the stopped sigma fields

(4.36) $$\mathcal{H}(s) \doteq \check{\mathcal{H}}(\tau(s)), \qquad s \geq 0.$$

Since $\{\tau(t) : t \geq 0\}$ is a nondecreasing sequence of stopping times, $\{\mathcal{H}(t) : t \geq 0\}$ is a filtration. Furthermore, the processes $\check{W}, \check{X}$ and $\check{U}$ are all linear combinations of the components of $\check{M}$ and $\check{U}$ and are hence adapted to $\{\check{\mathcal{H}}(t) : t \geq 0\}$. This shows that $W, X$ and $U$ are adapted to $\{\mathcal{H}(t) : t \geq 0\}$. We will now establish that $X$ is a $(0, \Sigma)$-Brownian motion with respect to this filtration, thereby completing the proof of Theorem 3.7.

Let $N = \mathbf{I} + \mathbf{J} + \mathbf{I} = 2\mathbf{I} + \mathbf{J}$ and $f \in \mathcal{C}_c^\infty(\mathbb{R}^N) =$ space of infinitely differentiable functions from $\mathbb{R}^N$ to $\mathbb{R}$ with compact support. Let $x = (x^E, x^S, x^\Phi) =$ $(x_i^E, i = 1, \ldots, \mathbf{I}; x_j^S, j = 1, \ldots, \mathbf{J}; x_i^\Phi, i = 1, \ldots, \mathbf{I}) \in \mathbb{R}^N$. Define the second order differential operator $G$ as follows:

$$Gf(x) = Gf(x^E, x^S, x^\Phi)$$

(4.37)
$$= \tfrac{1}{2} \sum_{i=1}^{\mathbf{I}} \sigma_{ii}^E f_{x_i^E x_i^E}(x) + \tfrac{1}{2} \sum_{j=1}^{\mathbf{J}} \sigma_{jj}^S f_{x_j^S x_j^S}(x)$$

$$+ \tfrac{1}{2} \sum_{i_1=1}^{\mathbf{I}} \sum_{i_2=1}^{\mathbf{I}} \sigma_{i_1 i_2}^\Phi f_{x_{i_1}^\Phi x_{i_2}^\Phi}(x),$$

where for $i = 1, \ldots, \mathbf{I}, \ j = 1, \ldots, \mathbf{J}, \ i_1, i_2 = 1, \ldots, \mathbf{I}$,

(4.38) $\quad \sigma_{ii}^E = \alpha_i (\alpha_i \sigma_i^u)^2, \qquad \sigma_{jj}^S = \beta_j x_j^* (\beta_j \sigma_j^v)^2, \qquad \sigma_{i_1 i_2}^\Phi = \sum_{j=1}^{\mathbf{J}} \beta_j x_j^* (\sigma_{i_1 i_2}^{\phi^j}).$

In order to complete the proof, it suffices to show that for any $t, s \geq 0$,

(4.39) $\quad \mathbb{E}\left[ f(M(t+s)) - f(M(t)) - \int_t^{t+s} Gf(M(u)) \, du \Big| \mathcal{H}(t) \right] = 0.$

We claim that, in order to prove (4.39), it is enough to prove the following:

(4.40)
$$\mathbb{E}\left[ h(\check{M}(s_k), \check{U}(s_k) : s_k \leq t, k = 1, \ldots, n) \Big( f(\check{M}(t+s)) - f(\check{M}(t)) - \int_t^{t+s} Gf(\check{M}(u)) \, d\tilde{\tau}(u) \Big) \right] = 0,$$

where $h$ is a bounded continuous function of its arguments. In order to justify the claim, define $\check{Y}_f(t) \doteq f(\check{M}(t)) - \int_0^t Gf(\check{M}(u)) \, d\tilde{\tau}(u), \ t \geq 0$. Then (4.40) implies that $\{\check{Y}_f(t)\}_{t \geq 0}$ is a $\{\check{\mathcal{H}}'(t)\}_{t \geq 0}$ martingale and hence a $\{\check{\mathcal{H}}(t)\}_{t \geq 0}$ martingale. We will now use optional sampling theorem to show that $\{\check{Y}_f(\tau(t))\}_{t \geq 0}$ is a $\{\mathcal{H}(t)\}_{t \geq 0}$ martingale. Using the fact that both $f$ and $Gf$ are bounded (by some constant $c > 0$) and (3.37), we have, for all $t \geq 0$,

$$\mathbb{E}[|\check{Y}_f(\tau(t))|]$$
$$\leq \mathbb{E}\left[ \left| f(\check{M}(\tau(t))) - \int_{[0, \tau(t)]} Gf(\check{M}(u)) \, d\tilde{\tau}(u) \right| \right]$$
$$\leq c(1+t),$$

(4.41) $\quad \mathbb{E}[|\check{Y}_f(T)| \mathbf{I}_{\{\tau(t) \geq T\}}]$
$$\leq \mathbb{E}\left[ \left| f(\check{M}(T)) - \int_{[0,T]} Gf(\check{M}(u)) \, d\tilde{\tau}(u) \right| \mathbf{I}_{\{\tilde{\tau}(T) \leq t\}} \right]$$
$$\leq c\mathbb{E}((1+t) \mathbf{I}_{\{\tilde{\tau}(T) \leq t\}})$$
$$\leq c(1+t) \mathbb{P}[\tilde{\tau}(T) \leq t] \to 0 \qquad \text{as } T \to \infty,$$

where the last convergence is a consequence of Lemma 3.6(B). Hence, it follows from Theorem 2.2.13 of [10], on observing that $\tau(t)$ is an a.s. finite $\{\tilde{\mathcal{H}}(s) : s \geq 0\}$ stopping time for each $t \geq 0$, that $\{\tilde{Y}_f(\tau(t))\}_{t \geq 0}$ is a $\{\mathcal{H}(t)\}_{t \geq 0}$ martingale. Thus, using Lemma 3.10, (4.39) follows. This proves the claim.

We now show (4.40). Fix $t, s \geq 0$. Let $t_m^r = t + \frac{m}{r^2}s$, $m = 0, 1, \ldots, r^2 - 1$. Define

$$\Psi^r \doteq \frac{1}{2} \sum_{i=1}^{\mathbf{I}} \sum_{m=0}^{r^2-1} f_{x_i^E x_i^E}(\check{N}^r(t_m^r))(\alpha_i^r \sigma_i^{u,r})^2 (\bar{E}_i^r(\check{\tau}^r(t_{m+1}^r)) - \bar{E}_i^r(\check{\tau}^r(t_m^r)))$$

$$+ \frac{1}{2} \sum_{j=1}^{\mathbf{J}} \sum_{m=0}^{r^2-1} f_{x_j^S x_j^S}(\check{N}^r(t_m^r))(\beta_j^r \sigma_j^{v,r})^2$$

(4.42)
$$\times (\bar{S}_j^r(\bar{T}_j^r(\check{\tau}^r(t_{m+1}^r))) - \bar{S}_j^r(\bar{T}_j^r(\check{\tau}^r(t_m^r))))$$

$$+ \frac{1}{2} \sum_{i_1=1}^{\mathbf{I}} \sum_{i_2=1}^{\mathbf{I}} \sum_{m=0}^{r^2-1} f_{x_{i_1}^\Phi x_{i_2}^\Phi}(\check{N}^r(t_m^r))$$

$$\times \sum_{j=1}^{\mathbf{J}} (\sigma_{i_1 i_2}^{\phi^j})(\bar{S}_j^r(\bar{T}_j^r(\check{\tau}^r(t_{m+1}^r))) - \bar{S}_j^r(\bar{T}_j^r(\check{\tau}^r(t_m^r)))).$$

Fix the subsequence along which the convergence to the limit point in Lemma 3.6(B) takes place and, by relabeling, denote it again by $\{r\}$. In view of (4.27), we have, in particular, that as $r \to \infty$,

(4.43) $\qquad (\check{M}^r, \check{N}^r, \check{\tau}^r, \check{U}^r) \Rightarrow (\check{M}, \check{M}, \check{\tau}, \check{U}).$

By appealing to the Skorohod representation theorem, assume without loss of generality that the above convergence holds a.s.

We will now show the following two results:

(4.44) $\qquad \limsup_{r \to \infty} \mathbb{E}\left[\left|\Psi^r - \int_t^{t+s} Gf(\check{M}(u))\, d\check{\tau}(u)\right|\right] = 0,$

(4.45)
$$\limsup_{r \to \infty} \mathbb{E}[h(\check{N}^r(s_k), \check{U}^r(s_k) : s_k \leq t, k = 1, \ldots, n)$$
$$\times (f(\check{N}^r(t+s)) - f(\check{N}^r(t)) - \Psi^r)] = 0.$$

Note that since both $h$ and $f$ are bounded and continuous, on applying the continuous mapping theorem and dominated convergence theorem, one gets (4.40) from (4.43) and (4.44)–(4.45). Thus the proof of the theorem will be complete once (4.44) and (4.45) are established.

We first consider (4.44). For $u \in [0, s)$, define

$$\check{N}_t^{r,*}(u) \doteq \sum_{m=0}^{r^2-1} \check{N}^r(t_m^r) \mathbf{I}_{[t_m^r, t_{m+1}^r)}(t + u).$$

From (4.43), it follows that as $r \to \infty$, $\check{N}_t^{r,*}(\cdot) \to \check{M}(t+\cdot)$ in $\mathcal{D}([0,s))$, a.s. Also, combining (4.43) with Lemma 3.6(C), we see that $\bar{E}_i^r(\check{\tau}^r(\cdot)) \to \alpha_i \check{\tau}(\cdot)$ in $\mathcal{D}$. Hence, using Lemma 3.11, we have the following result for the first term on the right side of (4.42):

$$(4.46) \quad \tfrac{1}{2} \sum_{i=1}^{\mathbf{I}} \sum_{m=0}^{r^2-1} f_{x_i^E x_i^E}(\check{N}^r(t_m^r))(\alpha_i^r \sigma_i^{u,r})^2 (\bar{E}_i^r(\check{\tau}^r(t_{m+1}^r)) - \bar{E}_i^r(\check{\tau}^r(t_m^r)))$$

$$= \tfrac{1}{2} \sum_{i=1}^{\mathbf{I}} \int_t^{t+s} f_{x_i^E x_i^E}(\check{N}_t^{r,*}(u-t))(\alpha_i^r \sigma_i^{u,r})^2 \, d\bar{E}_i^r(\check{\tau}^r(u))$$

$$\to \tfrac{1}{2} \sum_{i=1}^{\mathbf{I}} \int_t^{t+s} f_{x_i^E x_i^E}(\check{M}(u))(\alpha_i \sigma_i^u)^2 \alpha_i \, d\check{\tau}(u)$$

$$(4.47) \quad = \int_t^{t+s} \tfrac{1}{2} \sum_{i=1}^{\mathbf{I}} (\sigma_{ii}^E) f_{x_i^E x_i^E}(\check{M}(u)) \, d\check{\tau}(u),$$

almost surely. Using the bound $\mathbb{E}[(E_i^r(r^2 t)/r^2)^2] \leq c(t^2+1)$ for the renewal process $E_i^r$, where $c$ is a constant independent of $r$ and $t$, it follows that the expected value of the square of the term in (4.46) is bounded by

$$\kappa \mathbb{E}[(\bar{E}_i^r(\check{\tau}^r(t+s)))^2] \leq \kappa c((t+s)^2 + 1),$$

for some $\kappa > 0$. This estimate gives us the required uniform integrability enabling us to conclude that

$$(4.48) \quad \mathbb{E}\Bigg[\Bigg| \tfrac{1}{2} \sum_{i=1}^{\mathbf{I}} \sum_{m=0}^{r^2-1} f_{x_i^E x_i^E}(\check{N}^r(t_m^r))(\alpha_i^r \sigma_i^{u,r})^2 \\ \times (\bar{E}_i^r(\check{\tau}^r(t_{m+1}^r)) - \bar{E}_i^r(\check{\tau}^r(t_m^r))) \\ - \int_t^{t+s} \tfrac{1}{2} \sum_{i=1}^{\mathbf{I}} (\sigma_{ii}^E) f_{x_i^E x_i^E}(\check{M}(u)) \, d\check{\tau}(u) \, du \Bigg|\Bigg] \to 0.$$

Using similar arguments for the other two terms in $\Psi^r$ in (4.42) and the definition of the operator $G$ in (4.1), one gets (4.44).

Finally, we consider (4.45). Recall that $\check{U}^r(\cdot)$ and $\check{N}^r(\cdot)$ are adapted to $\{\mathcal{G}^r(t)\}$ [see Lemma 4.1 and remarks below (4.13)]. Thus, $h(\check{N}^r(s_k), \check{U}^r(s_k) : s_k \leq t, k = 1, \ldots, n)$ is measurable w.r.t. $\mathcal{G}^r(t)$. Hence, by boundedness of $h$ and $f$ and uniform integrability of $\Psi^r$ established above, in order to prove (4.45), it is enough to establish that

$$(4.49) \quad \limsup_{r \to \infty} |\mathbb{E}[f(\check{N}^r(t+s)) - f(\check{N}^r(t)) - \Psi^r | \mathcal{G}^r(t)]| = 0 \qquad \text{a.s.}$$

To prove (4.49), we first rewrite the difference $f(\check{N}^r(t+s)) - f(\check{N}^r(t))$ as a telescoping series and expand each term using Taylor's formula:

$$
\begin{aligned}
&\mathbb{E}[f(\check{N}^r(t+s)) - f(\check{N}^r(t)) - \Psi^r | \mathcal{G}^r(t)] \\
&= \mathbb{E}\left[\sum_{m=0}^{r^2-1} (f(\check{N}^r(t_{m+1}^r)) - f(\check{N}^r(t_m^r))) - \Psi^r \Big| \mathcal{G}^r(t)\right] \\
&\qquad(4.50) \\
&= \sum_{l,l'=1}^{N} \mathbb{E}\left[\sum_{m=0}^{r^2-1} (f_{x_l x_{l'}}(\theta_m^r) - f_{x_l x_{l'}}(\check{N}^r(t_m^r))) \right. \\
&\qquad\qquad \left. \times (\check{N}_l^r(t_{m+1}^r) - \check{N}_l^r(t_m^r))(\check{N}_{l'}^r(t_{m+1}^r) - \check{N}_{l'}^r(t_m^r)) \Big| \mathcal{G}^r(t)\right],
\end{aligned}
$$

where $\check{N}^r(\cdot) = (\check{N}_l^r(\cdot) : l = 1, \ldots, N)$ and $\theta_m^r$ lies in the line segment joining $\check{N}^r(t_m^r)$ and $\check{N}^r(t_{m+1}^r)$. In obtaining the last equality, we have used the fact that $\{\check{N}^r(t)\}_{t \geq 0}$ is a $\{\mathcal{G}^r(t)\}_{t \geq 0}$ martingale [see below (4.13)], and also applied (4.16)–(4.19). Define, for $c > 0$,

$$(4.51) \qquad \mathcal{B}_c^{r,m} \doteq \{|r(\check{N}^r(t_{m+1}^r) - \check{N}^r(t_m^r))| > c\}.$$

We claim that for each $\varepsilon > 0$, there exists a $c > 0$ such that

$$(4.52) \quad \sup_r \sup_{m=0,\ldots,r^2-1} \mathbb{E}[|r(\check{N}^r(t_{m+1}^r) - \check{N}^r(t_m^r))|^2 \mathbf{I}_{\mathcal{B}_c^{r,m}} | \mathcal{G}^r(t)] < \varepsilon.$$

Postponing the proof of (4.52), we now complete the proof of (4.49), assuming (4.52) holds. Note that for $\theta_m^r$ as in (4.50), we have

$$(4.53) \qquad |\theta_m^r - \check{N}^r(t_m^r)| < |\check{N}^r(t_{m+1}^r) - \check{N}^r(t_m^r)|.$$

Also, clearly, for all $l, l' = 1, \ldots, N$,

$$
\begin{aligned}
(4.54) \quad &|(\check{N}_l^r(t_{m+1}^r) - \check{N}_l^r(t_m^r))(\check{N}_{l'}^r(t_{m+1}^r) - \check{N}_{l'}^r(t_m^r))| \\
&\qquad < |\check{N}^r(t_{m+1}^r) - \check{N}^r(t_m^r)|^2.
\end{aligned}
$$

Since $f_{x_l x_{l'}}(\cdot)$ is bounded (by $c_{l,l'}$, say), we have, using (4.52), (4.53) and (4.54) in (4.50),

$$
\begin{aligned}
&|\mathbb{E}[(f(\check{N}^r(t+s)) - f(\check{N}^r(t)) - \Psi^r) | \mathcal{G}^r(t)]| \\
&\qquad \leq 2 \sum_{l,l'=1}^{N} \sum_{m=0}^{r^2-1} c_{l,l'} \mathbb{E}[|\check{N}^r(t_{m+1}^r) - \check{N}^r(t_m^r)|^2 \mathbf{I}_{\mathcal{B}_c^{r,m}} | \mathcal{G}^r(t)] \\
(4.55) &\qquad + \sum_{l,l'=1}^{N} \sum_{m=0}^{r^2-1} \left[\sup_{|x-y| < c/r} |f_{x_l x_{l'}}(x) - f_{x_l x_{l'}}(y)|\right]
\end{aligned}
$$

$$\times \mathbb{E}[|\check{N}^r(t^r_{m+1}) - \check{N}^r(t^r_m)|^2 \mathbf{I}_{(\mathcal{B}^{r,m}_c)^c} | \mathcal{G}^r(t)]$$

$$\leq 2\varepsilon \left( \sum_{l,l'=1}^{N} c_{l,l'} \right) + \sum_{l,l'=1}^{N} \sum_{m=0}^{r^2-1} \left[ \sup_{|x-y|<c/r} |f_{x_l x_{l'}}(x) - f_{x_l x_{l'}}(y)| \right] \frac{c^2}{r^2}.$$

Taking lim sup as $r \to \infty$ in (4.55), recalling that $f \in \mathcal{C}_c^\infty(\mathbb{R}^N)$ and that $\varepsilon > 0$ is arbitrary, the statement in (4.49) follows.

Finally, we prove the claim made in (4.52). Note that it is enough to prove that the sequence $\{|r(\check{N}^r_l(t^r_{m+1}) - \check{N}^r_l(t^r_m))|^2\}$ is uniformly integrable [conditioned on $\mathcal{G}^r(t)$] for each $l = 1, \ldots, N$. We will only check this for $i = 1, \ldots, \mathbf{I}$ (i.e., $l = 1, \ldots, \mathbf{I}$) here: the proof for other values of $l$ (i.e., $l = \mathbf{I} + 1, \ldots, 2\mathbf{I} + \mathbf{J} = N$) is similar and is omitted. Thus, we will show that for $i = 1, \ldots, \mathbf{I}$, given any $\varepsilon > 0$, there exists $c > 0$ such that

$$\sup_r \sup_{m=0,\ldots,r^2-1} \mathbb{E}[|r(\check{N}^{E,r}_i(t^r_{m+1}) - \check{N}^{E,r}_i(t^r_m))|^2$$
(4.56)
$$\times \mathbf{I}_{\{|r(\check{N}^{E,r}_i(t^r_{m+1}) - \check{N}^{E,r}_i(t^r_m))|>c\}} | \mathcal{G}^r(t)] < \varepsilon \quad \text{a.s.}$$

From the definitions of $\check{N}^{E,r}_i(\cdot)$ and $\tilde{\xi}^r_i(\cdot)$ [see (4.13) and (4.11)], we get that for $i = 1, \ldots, \mathbf{I}'$ [the case $i = \mathbf{I}' + 1, \ldots, \mathbf{I}$, (4.56) holds trivially],

(4.57) $$\begin{aligned} & r(\check{N}^{E,r}_i(t^r_{m+1}) - \check{N}^{E,r}_i(t^r_m)) \\ &= r(\tilde{\xi}^r_i(\rho^{r,E}_i(t^r_{m+1})) - \tilde{\xi}^r_i(\rho^{r,E}_i(t^r_m))) \\ &= (\rho^{r,E}_i(t^r_{m+1}) - \rho^{r,E}_i(t^r_m)) - \alpha^r_i(\xi^r_i(\rho^{r,E}_i(t^r_{m+1})) - \xi^r_i(\rho^{r,E}_i(t^r_m))). \end{aligned}$$

Using the inequality $\check{\tau}^r(t+s) - \check{\tau}^r(t) \leq s$ for all $s, t \geq 0$ and recalling the relation between $E^r_i$ and $\xi^r_i$ [see (2.5)], it follows that

(4.58) $$\begin{aligned} |\rho^{r,E}_i(t^r_{m+1}) - \rho^{r,E}_i(t^r_m)| &\leq E^r_i(r^2 \check{\tau}^r(t^r_m) + s) - E^r_i(r^2 \check{\tau}^r(t^r_m)) \\ &\leq E^r_i(\xi^r_i(\rho^{r,E}_i(t^r_m)) + s) - E^r_i(\xi^r_i(\rho^{r,E}_i(t^r_m))) + 1. \end{aligned}$$

Observing that $\rho^{r,E}_i(t^r_m)$ is a $\{\mathcal{F}^r((m,n))\}$ stopping time and using properties of renewal processes, we obtain

(4.59) $$\begin{aligned} &\mathbb{P}[E^r_i(\xi^r_i(\rho^{r,E}_i(t^r_m)) + s) - E^r_i(\xi^r_i(\rho^{r,E}_i(t^r_m))) \in \cdot | \mathcal{G}^r(t)] \\ &\stackrel{d}{<} \mathbb{P}[E^r_i(s) + 1 \in \cdot], \end{aligned}$$

where for probability measures $\mu, \nu$ on $\mathbb{R}$, we write $\mu \stackrel{d}{<} \nu$ if $\mu(-\infty, x] \geq \nu(-\infty, x]$ for all $x \in \mathbb{R}$. Proof of (4.59), which is left to the reader, can be given along the lines of [22] (page 63). This, in particular, implies the

CONTROLLED STOCHASTIC NETWORKS 41

following moment estimate for the term on the left-hand side of (4.58):

$$\mathbb{E}((\rho_i^{r,E}(t_{m+1}^r) - \rho_i^{r,E}(t_m^r))^2 \mathbf{I}_{\{|\rho_i^{r,E}(t_{m+1}^r) - \rho_i^{r,E}(t_m^r)| > c\}} | \mathcal{G}^r(t))$$
(4.60)
$$\leq \mathbb{E}((E_i^r(s) + 2)^2 \mathbf{I}_{\{E_i^r(s) + 2 > c\}}).$$

Using Lemma 4.2, we now have that the supremum (over $r$) of the right-hand side of the above inequality converges to zero as $c \to \infty$. This proves the (conditional) uniform integrability of the sequence corresponding to the square of the first term in (4.57). Thus, it suffices to prove that the square of the second term in (4.57) is (conditionally) uniformly integrable. For this, note that if $\rho_i^{r,E}(t_{m+1}^r) > \rho_i^{r,E}(t_m^r)$, then using the fact that $\check{\tau}^r(t_{m+1}^r) - \check{\tau}^r(t_m^r) \leq (t + (m+s)/r^2) - (t + ms/r^2) \leq s/r^2$, it follows that

$$\xi_i(\rho_i^{r,E}(t_{m+1}^r)) - \xi_i(\rho_i^{r,E}(t_m^r)) \leq r^2(\check{\tau}^r(t_{m+1}^r) - \check{\tau}^r(t_m^r)) + u_i^r(\rho_i^{r,E}(t_{m+1}^r))$$
(4.61)
$$\leq s + u_i^r(\rho_i^{r,E}(t_{m+1}^r)).$$

For the remaining case, that is, if $\rho_i^{r,E}(t_{m+1}^r) = \rho_i^{r,E}(t_m^r)$, the difference between the two terms in the left-hand side of (4.61) is zero. Hence, we have

$$[\xi_i(\rho_i^{r,E}(t_{m+1}^r)) - \xi_i(\rho_i^{r,E}(t_m^r))]^2 \mathbf{I}_{\{|\xi_i(\rho_i^{r,E}(t_{m+1}^r)) - \xi_i(\rho_i^{r,E}(t_m^r))| > c\}}$$
(4.62)
$$\leq 2[s^2 + (u_i^r(\rho_i^{r,E}(t_{m+1}^r)))^2] \mathbf{I}_{\{u_i^r(\rho_i^{r,E}(t_{m+1}^r)) > c - s\}}$$
$$\leq 2[\xi_i^c(\rho_i^{r,E}(t_{m+1}^r)) - \xi_i^c(\rho_i^{r,E}(t_m^r))],$$

where $\xi_i^c(m) \doteq \sum_{k=1}^{m} [s^2 + u_i^r(k)^2] \mathbf{I}_{\{u_i^r(k) > c - s\}}$ for $m \geq 1$. Note that $\{[s^2 + u_i^r(k)^2] \mathbf{I}_{\{u_i^r(k) > c - s\}}\}_{k \geq 1}$ is a sequence of i.i.d. random variables with finite mean and $\xi_i^c((m, n)) \equiv \xi_i^c(m_i)$ is an $\mathcal{F}^r((m, n))$ martingale. Since $\rho_i^{r,E}(t_{m+1}^r) \geq \rho_i^{r,E}(t_m^r)$ are two $\mathcal{F}^r((m, n))$ stopping times, using the optional sampling theorem, one gets the following from (4.62),

$$\mathbb{E}[\xi_i(\rho_i^{r,E}(t_{m+1}^r)) - \xi_i(\rho_i^{r,E}(t_m^r))]^2 \mathbf{I}_{\{|\xi_i(\rho_i^{r,E}(t_{m+1}^r)) - \xi_i(\rho_i^{r,E}(t_m^r))| > c\}} | \mathcal{G}^r(t_m^r)]$$
(4.63)
$$\leq 2\mathbb{E}(\rho_i^{r,E}(t_{m+1}^r) - \rho_i^{r,E}(t_m^r) | \mathcal{G}^r(t_m^r)) \mathbb{E}[(s^2 + u_i^r(1)^2) \mathbf{I}_{\{u_i^r(1) > c - s\}}]$$
$$\leq 2\mathbb{E}(E_i^r(s) + 2) \mathbb{E}[(s^2 + u_i^r(1)^2) \mathbf{I}_{\{u_i^r(1) > c - s\}}],$$

where the second inequality follows on using (4.58) and (4.59). Finally, since $\{u_i^r(1)^2\}_r$ is uniformly integrable [see (2.2)], and the expected value of $E_i^r(s)$ is bounded in $r$ (see Lemma 4.2), the right-hand side above goes to zero as $c \to \infty$. Thus, the square of the second term in (4.57) is (conditionally) uniformly integrable. This proves (4.56), and hence the proof of Theorem 3.7 is complete. □

**5. Admissibility of control policies.** In this section, we will show that for a very broad class of sequences of control policies $\{T^r\}$ which satisfy some natural nonanticipativity conditions, the admissibility requirement (iv) of Definition 2.6 is satisfied. Throughout this section, we will fix $r$ and consider the $r$th network $\mathcal{N}^r$. Consider a control policy $T^r$ which satisfies (i), (ii) and (iii) of Definition 2.6. We will show that under further natural conditions on $T^r$, the property (iv) of Definition 2.6 also holds. For the rest of the section, we will suppress $r$ from the notation, unless it is necessary. The first condition on the policy $T$ is the following:

ASSUMPTION 5.1. There is a **J**-dimensional measurable process $\dot{T}$ with values in $\{0,1\}$ such that $T(t) = \int_0^t \dot{T}(s)\,ds$ for all $t \in [0,\infty)$.

The second condition on the policy states that the process $\dot{T}$ does not change values between two successive event times, where an event is either an exogenous arrival into the system or the completion of a service by some server. To state this condition precisely, define $\Upsilon_0 = 0$, and for any $\ell \geq 0$, let $\Upsilon_{\ell+1}$ be the first time after $\Upsilon_\ell$ when either an arrival or a service completion takes place. From (2.1), it follows that, almost surely, $\{\Upsilon_\ell\}_{\ell \in \mathbb{N}_0}$ is a strictly increasing sequence, increasing to $\infty$ as $\ell \to \infty$. Here we assume that multiple events can occur at a given $\Upsilon_\ell$.

ASSUMPTION 5.2. For all $\ell \in \mathbb{N}_0$, $\dot{T}(t) = \dot{T}(\Upsilon_\ell)$ for all $t \in [\Upsilon_\ell, \Upsilon_{\ell+1})$.

Our final condition on the policy is a natural nonanticipativity property. In order to state this condition, we introduce the following notation. For $i \in \{1,\ldots,\mathbf{I}\}$ and $\ell \in \mathbb{N}_0$, let $u_i^\ell \doteq \xi_i(E_i(\Upsilon_\ell)+1) - \Upsilon_\ell$. Thus, $u_i^\ell$ is the residual (exogenous) arrival time at the $i$th buffer at time $\Upsilon_\ell$, unless an arrival of the $i$th class occurred at time $\Upsilon_\ell$, in which case it equals $u_i^\ell = u_i(E_i(\Upsilon_\ell)+1)$. Similarly, for $j \in \{1,\ldots,\mathbf{J}\}$, $\ell \in \mathbb{N}_0$, define $v_j^\ell \doteq \eta_j(S_j(\Upsilon_\ell)+1) - \Upsilon_\ell$. Next, for $i \in \{1,\ldots,\mathbf{I}\}$, set $Q_{i,0} = 0$, and for $\ell \geq 1$, $Q_{i,\ell} \doteq Q_i(\Upsilon_\ell)$. Also, for $j \in \{1,\ldots,\mathbf{J}\}$ and $\ell \geq 0$, let $\dot{T}_j^\ell \doteq \dot{T}_j(\Upsilon_\ell)$. Also, let $\dot{T}_j^{-1} \doteq 0$. Finally, define, for $\ell \geq 0$,

$$(5.1) \quad \chi^\ell \doteq \{(\Upsilon_{\ell'}, u_i^{\ell'}, v_j^{\ell'}, Q_i^{\ell'}, \dot{T}_j^{\ell'-1} : i \in 1,\ldots,\mathbf{I}, j \in 1,\ldots,\mathbf{J}) : \ell' = 0,\ldots,\ell\}.$$

ASSUMPTION 5.3. For all $\ell \geq 0$, $\dot{T}(\Upsilon_\ell)$ is a measurable function of $\chi^\ell$.

The following is the main result of this section:

THEOREM 5.4. *Let $T$ be a control policy for $\mathcal{N}^r$ that, in addition to* (i), (ii) *and* (iii) *of Definition* 2.6, *satisfies Assumptions* 5.1, 5.2 *and* 5.3. *Then $T$ satisfies condition* (iv) *of Definition* 2.6.

PROOF. We begin with some notation. For a given activity $j \in \{1, \ldots, \mathbf{J}\}$, let $\iota(j)$ denote the associated buffer. For a given buffer $i \in \{1, \ldots, \mathbf{I}\}$, let $\mathcal{J}_i$ denote the set of all activities that are associated with the buffer $i$. Define, for each $\tilde{\mathcal{J}}_i \subseteq \mathcal{J}_i$ and a vector $a \in \{0,1\}^{\mathbf{J}}$,

$$(5.2) \qquad \pi(a, \tilde{\mathcal{J}}_i) \doteq \left(\prod_{j \in \tilde{\mathcal{J}}_i} a_j\right)\left(\prod_{j \in \mathcal{J}_i \setminus \tilde{\mathcal{J}}_i}(1-a_j)\right).$$

Note that for $\tilde{\mathcal{J}}_i \subseteq \mathcal{J}_i$, $\pi(a, \tilde{\mathcal{J}}_i)$ is either zero or one, and it takes the value one only if $a_j = 1$ for all $j \in \tilde{\mathcal{J}}_i$ and $a_j = 0$ for all $j \in \mathcal{J}_i \setminus \tilde{\mathcal{J}}_i$. Denote this unique $\tilde{\mathcal{J}}_i$ by $\mathcal{J}_i(a)$. A straightforward calculation shows that

$$(5.3) \qquad \Upsilon_{\ell+1} = \Upsilon_\ell + \min_{i \in \{1, \ldots, \mathbf{I}\}} \min\{u_i^\ell, v_j^\ell : j \in \mathcal{J}_i(\dot{T}^\ell)\}.$$

Also for $i = 1, \ldots, \mathbf{I}$, let $\mathcal{I}_i^{\ell+1}$ be the indicator function that at $\Upsilon_{\ell+1}$, the event (arrival or service completion) occurred at buffer $i$. More precisely, for $i = 1, \ldots, \mathbf{I}$ and $\ell \geq 0$,

$$(5.4) \qquad \mathcal{I}_i^{\ell+1} = \begin{cases} 1, & \text{if } \min\{u_i^\ell, v_j^\ell : j \in \mathcal{J}_i(\dot{T}^\ell)\} \\ & \quad = \min_{i' \in \{1, \ldots, \mathbf{I}\}} \min\{u_{i'}^\ell, v_j^\ell : j \in \mathcal{J}_{i'}(\dot{T}^\ell)\}, \\ 0, & \text{otherwise.} \end{cases}$$

From (5.4), (5.3) and Assumption 5.3, it follows that both

$$(5.5) \qquad \mathcal{I}_i^{\ell+1}, \Upsilon_{\ell+1} \text{ are measurable functions of } \chi^\ell.$$

For $\ell \geq 0$ and $(m, n) \in \mathbb{N}^{\mathbf{I}+\mathbf{J}}$, let

$$(5.6) \qquad \begin{aligned} \mathcal{B}_{m,n}^\ell &\doteq \{E_i(\Upsilon_\ell) + 1 = m_i, \\ & S_j(T_j(\Upsilon_\ell)) + 1 = n_j : i = 1, \ldots, \mathbf{I}, \ j = 1, \ldots, \mathbf{J}\}. \end{aligned}$$

We claim that the following two statements hold for all $\ell = 0, 1, \ldots$:

$$(5.7) \qquad \begin{aligned} &\text{(A)} \quad \mathcal{B}_{m,n}^\ell \in \mathcal{F}((m,n)) \text{ for all } (m,n) \in \mathbb{N}^{\mathbf{I}+\mathbf{J}}, \\ &\text{(B)} \quad \mathbf{I}_{\mathcal{B}_{m,n}^\ell} \chi^\ell \in \mathcal{F}((m,n)) \text{ for all } (m,n) \in \mathbb{N}^{\mathbf{I}+\mathbf{J}}, \end{aligned}$$

where $\mathbf{I}_C$ denotes the indicator function of a set $C$. We will show (5.7) by induction on $\ell$. Note that for $\ell = 0$, $\mathcal{B}_{m,n}^\ell$ is the sample space if $(m,n) = \mathbf{1}$, otherwise it is the empty set. Thus, for $\ell = 0$, (A) holds trivially. Also, when $(m,n) = \mathbf{1}$, one sees from (5.1) that $\chi^0 \in \mathcal{F}((m,n))$. So (B) also holds for the case $\ell = 0$. Now suppose that (5.7) holds for some $\ell \in \mathbb{N}_0$, and consider (5.7) with $\ell$ replaced by $\ell + 1$.

Note that the set $\mathcal{B}_{m,n}^{\ell+1}$ can be represented as

$$(5.8) \qquad \mathcal{B}_{m,n}^{\ell+1} = \bigcup_{(p,q) \leq (m,n)} \mathcal{B}_{p,q}^\ell \cap \mathcal{B}_{m,n}^{\ell+1}.$$

We will next show that

$$(5.9) \qquad \mathcal{B}_{p,q}^{\ell} \cap \mathcal{B}_{m,n}^{\ell+1} \in \mathcal{F}((m,n)) \qquad \text{for all } p \leq m, \ q \leq n.$$

This, in view of (5.8), will prove that (A) in (5.7) holds with $\ell$ replaced by $\ell+1$. By induction, $\mathcal{B}_{p,q}^{\ell} \in \mathcal{F}((p,q)) \subseteq \mathcal{F}((m,n))$ for all $p \leq m, \ q \leq n$. Also, on $\mathcal{B}_{p,q}^{\ell}$, the following representations hold. For $i \in \{1,\ldots,\mathbf{I}\}$, $j \in \{1,\ldots,\mathbf{J}\}$,

$$(5.10) \qquad [E_i(\Upsilon_{\ell+1})+1] = p_i + \mathcal{I}_i^{\ell+1}\mathbf{I}_{\{u_i^{\ell} < v_j^{\ell} : j \in \mathcal{J}_i(\dot{T}^{\ell})\}},$$

$$(5.11) \ \ [S_j(T_j(\Upsilon_{\ell+1}))+1] = q_j + \dot{T}_j^{\ell}\mathcal{I}_{\imath(j)}^{\ell+1}\mathbf{I}_{\{v_j^{\ell} < u_{\imath(j)}^{\ell}, v_{j'}^{\ell} : j' \in \mathcal{J}_{\imath(j)}(\dot{T}^{\ell}), j' \neq j\}}.$$

Noting that the right-hand sides of (5.10) and (5.11) are measurable functionals of $\chi^{\ell}$, we see from (5.6) that

$$(5.12) \qquad \mathbf{I}_{\mathcal{B}_{p,q}^{\ell} \cap \mathcal{B}_{m,n}^{\ell+1}} = \mathbf{I}_{\mathcal{B}_{p,q}^{\ell}} \Psi_{\ell}(\chi^{\ell})$$

for some measurable function $\Psi_{\ell}$. Now by induction, the right-hand side of (5.12) is in $\mathcal{F}((p,q))$ and therefore in $\mathcal{F}((m,n))$. This proves (5.9) and hence (A) in (5.7) holds with $\ell$ replaced by $\ell+1$.

Finally, we verify (B) in (5.7) for $\ell+1$. Note that on $\mathcal{B}_{m,n}^{\ell+1}$, for $i = 1,\ldots,\mathbf{I}$,

$$(5.13) \quad Q_i^{\ell+1} = q^r + (m_i - 1) - (n_j - 1)\sum_{j=1}^{\mathbf{J}} C_{ij} + \sum_{j=1}^{\mathbf{J}} \mathbf{\Phi}_i^{j,r}(n_j - 1).$$

Thus, using (5.7)(A) for $\ell+1$, we see that $Q_i^{\ell+1}\mathbf{I}_{\mathcal{B}_{m,n}^{\ell+1}} \in \mathcal{F}((m,n))$. Recalling that $\Upsilon_{\ell+1}$ and $\dot{T}^{\ell}$ are measurable functions of $\chi^{\ell}$, we see from (5.8), the induction hypothesis and (5.7)(A) for $\ell+1$ that $\mathbf{I}_{\mathcal{B}_{m,n}^{\ell+1}}(\Upsilon_{\ell+1}, \dot{T}^{\ell}) \in \mathcal{F}((m,n))$.

Hence, we have shown

$$(5.14) \quad \mathbf{I}_{\mathcal{B}_{m,n}^{\ell+1}}(\Upsilon_{\ell+1}, Q_i^{\ell+1}, \dot{T}_j^{\ell} : i \in \{1,\ldots,\mathbf{I}\}, j \in \{1,\ldots,\mathbf{J}\}) \in \mathcal{F}((m,n)).$$

Next, on $\mathcal{B}_{m,n}^{\ell+1}$, for $i \in \{1,\ldots,\mathbf{I}\}, \ j \in \{1,\ldots,\mathbf{J}\}$, we have

$$(5.15) \ \ u_i^{\ell+1} = \mathbf{I}_{\{u_i^{\ell} > \Upsilon_{\ell+1} - \Upsilon_{\ell}\}}(u_i^{\ell} - \Upsilon_{\ell+1} + \Upsilon_{\ell}) + \mathbf{I}_{\{u_i^{\ell} = \Upsilon_{\ell+1} - \Upsilon_{\ell}\}}u_i(m_i),$$

$$v_j^{\ell+1} = (1-\dot{T}_j^{\ell})v_j^{\ell} + \dot{T}_j^{\ell}(1 - \mathcal{I}_{\imath(j)}^{\ell+1})(v_j^{\ell} - \Upsilon_{\ell+1} + \Upsilon_{\ell})$$

$$(5.16) \qquad\quad + \dot{T}_j^{\ell}\mathcal{I}_{\imath(j)}^{\ell+1}\mathbf{I}_{\{v_j^{\ell} \leq u_{\imath(j)}^{\ell}, v_{j'}^{\ell} : j' \in \mathcal{J}_{\imath(j)}(\dot{T}^{\ell}), j' \neq j\}}v_j(n_j)$$

$$\qquad\quad + \dot{T}_j^{\ell}\mathcal{I}_{\imath(j)}^{\ell+1}(1 - \mathbf{I}_{\{v_j^{\ell} \leq u_{\imath(j)}^{\ell}, v_{j'}^{\ell} : j' \in \mathcal{J}_{\imath(j)}(\dot{T}^{\ell}), j' \neq j\}})(v_j^{\ell} - \Upsilon_{\ell+1} + \Upsilon_{\ell}).$$

The above two identities, together with (5.8), (5.5) and the fact that (5.7)(A) holds with $\ell+1$, show that $\mathbf{I}_{\mathcal{B}_{m,n}^{\ell+1}}(u_i^{\ell+1}, v_j^{\ell+1}) \in \mathcal{F}((m,n))$. Combining this with (5.14) and using the induction hypothesis [along with (5.7)(A) for $\ell+1$], we now have that (5.7)(B) holds for $\ell+1$. This proves (5.7) for all $\ell \in \mathbb{N}_0$.

Finally, using (5.5) and (5.7), we have, for each fixed $t \geq 0$ and all $(m,n) \in \mathbb{N}^{\mathbf{I}+\mathbf{J}}$,

$$\begin{aligned}
\mathcal{B}_{m,n} &\doteq \{E_i(t) + 1 = m_i, S_j(T_j(t)) + 1 = n_j : i = 1, \ldots, \mathbf{I}, j = 1, \ldots, \mathbf{J}\} \\
&= \bigcup_{\ell=1}^{p} \mathcal{B}_{m,n}^{\ell} \cap \{\Upsilon_\ell \leq t < \Upsilon_{\ell+1}\} \in \mathcal{F}((m,n)),
\end{aligned} \tag{5.17}$$

where $p \doteq \sum_{i=1}^{\mathbf{I}} m_i + \sum_{j=1}^{\mathbf{J}} n_j$. This proves the first part of condition (iv) in Definition 2.6. To show the final part of this condition, it suffices to prove that for $0 \leq t_j \leq t$, $j = 1, \ldots, \mathbf{J}$, and all $(m,n) \in \mathbb{N}^{\mathbf{I}+\mathbf{J}}$,

$$\{T_j(t) < t_j, j = 1, \ldots, \mathbf{J}\} \cap \mathcal{B}_{m,n} \in \mathcal{F}((m,n)). \tag{5.18}$$

From the properties of $T$ [Assumptions (5.1), (5.2)], we have

$$\begin{aligned}
\{T_j(t) &< t_j, j = 1, \ldots, \mathbf{J}\} \cap \mathcal{B}_{m,n} \\
&= \bigcup_{\ell=1}^{p} \mathcal{B}_{m,n}^{\ell} \cap \{\Upsilon_\ell \leq t < \Upsilon_{\ell+1}\} \\
&\quad \cap \{T_j(\Upsilon_\ell) + \dot{T}_j(\Upsilon_\ell)(t - \Upsilon_\ell) < t_j, j = 1, \ldots, \mathbf{J}\} \\
&\in \mathcal{F}((m,n)),
\end{aligned} \tag{5.19}$$

where the last inclusion follows on using (5.7) and (5.5). This proves (5.18) and, hence, the theorem. $\square$

DEPARTMENT OF STATISTICS
AND OPERATIONS RESEARCH
UNIVERSITY OF NORTH CAROLINA
CHAPEL HILL, NORTH CAROLINA 27599-3260
USA
E-MAIL: budhiraj@email.unc.edu

DEPARTMENT OF STATISTICS
IOWA STATE UNIVERSITY
303 SNEDECOR HALL
AMES, IOWA 50011-1210
USA
E-MAIL: apghosh@iastate.edu